\documentclass[12pt]
{amsart}
\usepackage{latexsym,amscd,amssymb}  
\makeatletter
\theoremstyle{plain}
  \newtheorem{thm}{Theorem}[section]
  \newtheorem{prop}[thm]{Proposition}
  \newtheorem{lem}[thm]{Lemma}
  \newtheorem{cor}[thm]{Corollary}
\theoremstyle{definition}
  \newtheorem{dfn}[thm]{Definition}
  \newtheorem{exmp}[thm]{Example}
  
\theoremstyle{remark}
  \newtheorem{rem}[thm]{Remark}
\numberwithin{equation}{section}

\def\ba{{\bf a}}
\def\bb{{\bf b}}
\def\bc{{\bf c}}
\def\be{{\bf e}}
\def\ff{{\bf f}}
\def\bx{{\bf x}}
\def\00{{\bf 0}}
\def\b1{{\bf 1}}
\def\E{{}^*E}
\def\Dc{\mathcal{D}^\bullet}
\def\bP{ {\mathbf P}}
\def\bH{{\mathbb H}}
\def\P{P^\bullet}
\def\Cb{\overline{C}}
\def\R{{\mathbb R}}
\def\Q{{\mathbb Q}}
\def\Z{{\mathbb Z}}
\def\NN{{\mathbb N}}
\def\cL{{\mathbf L}}
\def\Dcom{\mathcal{D}^\bullet}
\def\Ddel{\Dcom_X}
\def\Dsig{\Dcom_Y}
\def\wdel{\omega^\bullet_{K[\Delta]}}
\def\Wde{\omega_{K[\Delta]}}
\def\cExt{{\mathcal Ext}}
\def\cHom{{\mathcal Hom}}

\def\I{I^\bullet}
\def\J{J^\bullet}
\def\M{M^\bullet}
\def\N{N^\bullet}
\def\T{T^\bullet}
\def\L{L^\bullet}
\def\Sh{\operatorname{Sh}}
\def\const{\underline{K}}
\def\cF{{\mathcal F}}
\def\cG{{\mathcal G}^\bullet}
\def\rH{\tilde{H}}
\def\cH{{\mathcal H}}
\def\m{{\mathfrak m}}
\def\chara{\operatorname{char}}
\def\grR{\operatorname{*mod}R}
\def\CgrR{\operatorname{*mod}_CR}
\def\GrR{\operatorname{*Mod}R}
\def\Ass{\operatorname{Ass}}
\def\ann{\operatorname{ann}}

\def\irrdim{\operatorname{ir-dim}}
\def\lin{\operatorname{lin}}
\def\nat{\operatorname{nat}} 
\def\Sq{\operatorname{Sq}}
\def\supp{\operatorname{supp}_+}
\def\Supp{\operatorname{supp}}
\def\inj*dim{\operatorname{*inj.dim}_R}
\def\<{{\langle}}
\def\>{{\rangle}}
\def\Ext{\operatorname{Ext}}
\def\Hom{\operatorname{Hom}}
\def\depth{\operatorname{depth}}
\def\lcm{\operatorname{lcm}}
\def\can{\omega_R}
\def\wS{\omega_S}
\def\RHom{{\rm R}{\Hom}}
\def\RcHom{{\rm R}{\cHom}}
\def\Du{{\mathbf D}}
\def\Id{\operatorname{Id}}
\def\relint{\operatorname{rel-int}}

\begin{document}

\title{Notes on $C$-graded modules over an affine semigroup ring $K[C]$}
\author{Kohji Yanagawa}
\address{Department of Mathematics, 
Graduate School of Science, Osaka University, Toyonaka, Osaka 
560-0043, Japan}
\email{yanagawa@math.sci.osaka-u.ac.jp}

\maketitle
\begin{abstract}
Let $C \subset \NN^d$ be an affine semigroup, and $R=K[C]$ its 
semigroup ring. This paper is a collection of various results on 
``$C$-graded" $R$-modules $M = \bigoplus_{\bc \in C} M_\bc$, 
especially, monomial ideals of $R$.
For example, we show the following: If $R$ is normal 
and $I \subset R$ is a {\it radical} monomial ideal 
 (i.e., $R/I$ is a generalization of Stanley-Reisner rings), then  
the sequentially Cohen-Macaulay property of $R/I$ 
is a topological property of the ``geometric realization" 
of the cell complex associated with $I$. Moreover, we can give a  
squarefree modules/constructible sheaves version of this result. 
We also show that if $R$ is normal and $I \subset R$ 
is a Cohen-Macaulay monomial ideal then $\sqrt{I}$ is Cohen-Macaulay again. 
\end{abstract}

\section{Introduction}
First, we fix the notation used throughout this paper. 
Let $C \subset \Z^d \subset \R^d$ be an affine semigroup 
(i.e., $C$ is a finitely generated additive submonoid of $\Z^d$), 
and $R := K[\bx^\bc \mid \bc \in C] \subset K[x_1^{\pm 1}, \ldots, 
x_d^{\pm 1}]$ the semigroup ring of $C$ over a field $K$. 
Here $\bx^\bc$ for $\bc = (c_1, \ldots, c_d) \in C$  
denotes the monomial $\prod_{i=1}^dx_i^{c_i}$.  
We always assume that $\Z C = \Z^d$ and $C \cap (-C) = \{ 0 \}$. 
Thus $\dim R = d$ and 
$\m := (\bx^\bc \mid 0 \ne \bc \in C)$ is a maximal ideal of $R$. 

Of course, $R = \bigoplus_{\bc \in C} K \bx^\bc$ is a $\Z^d$-graded ring. 
We say a $\Z^d$-graded ideal of $R$ is a {\it monomial ideal}. 
Let $\GrR$ be the category of $\Z^d$-graded $R$-modules and their degree 
preserving $R$-homomorphisms, and $\grR$ its full subcategory consisting 
of finitely generated modules. As usual, for $M \in \GrR$ and 
$\ba \in \Z^d$, $M_\ba$ denotes the degree $\ba$ component of $M$, and 
$M(\ba)$ denotes the 
shifted module of $M$ with $M(\ba)_\bb = M_{\ba + \bb}$. 
We say $M \in \GrR$ is {\it $C$-graded}, 
if $M_\ba = 0$ for all $\ba \not \in C$. A monomial ideal $I \subset R$ 
and the quotient ring $R/I$ are $C$-graded modules. 
Let $\CgrR$ be the full subcategory of $\grR$ 
consisting of $C$-graded modules.  

Miller~\cite{Mil} proved that $\CgrR$ has enough injectives and any object 
has a minimal injective resolution in this category, which is unique up to 
isomorphism and has finite length. This resolution is called a  
{\it minimal irreducible resolution}, since an indecomposable injective in 
$\CgrR$ corresponds to a monomial irreducible ideal. 

In \S2, under the assumption that $R$ is Cohen-Macaulay and simplicial,  
we show that information on $M \in \CgrR$ such as depth and  
Cohen-Macaulay property can be read off 
from numerical invariants of the minimal irreducible resolution of $M$ 
 (something analogous to ``Bass numbers"). 
One might think these results should be a variant of the fact that 
depth and related conditions can be characterized by Bass numbers. 
But this insight is not quite correct. 
Philosophically, our result is rather closer to a 
theorem of Eagon and Reiner~\cite{ER} stating that the Stanley-Reisner ring 
of a simplicial complex $\Delta$ is Cohen-Macaulay if and only if 
that of the  {\it Alexander dual} $\Delta^*$ 
has a linear free resolution, and Miller's generalization of this result 
to finitely generated $\NN^d$-graded modules over 
a polynomial ring $K[x_1, \ldots, x_d]$ (see \cite{Mil0}). 
In fact, in the polynomial ring case, the ``Alexander dual" of 
our Theorem~\ref{depth} corresponds to his \cite[Theorem~4.20]{Mil0}. 
But the proofs are not similar. 

In \S\S3-6, we assume that $R$ is normal (but not necessarily simplicial).  
In \S3, we study the full subcategory $\Sq$ of $\CgrR$ consisting of 
{\it squarefree modules}. 
The notion of squarefree modules over a normal semigroup ring was 
introduced by the author in \cite{Y2}.  
A monomial ideal $I \subset R$ is squarefree 
if and only if it is a radical ideal (i.e., $I = \sqrt{I}$). 
The category $\Sq$ behaves much nicer than $\CgrR$ 
as shown in \cite{Y2}. In this paper, we show that a squarefree 
module $M$ is  {\it sequentially Cohen-Macaulay} 
(a non-pure generalization of the Cohen-Macaulay property) 
if and only if the 
{\it linear strand} of its minimal irreducible resolution is acyclic. 
 
In \S4, assuming that $R$ is normal,  
we study the quotient ring $R/I$ by a {\it radical} monomial ideal 
$I$, focusing on the problem when $R/I$ is  
sequentially Cohen-Macaulay. 
When $R$ is a polynomial ring, $R/I$ is usually called a 
{\it Stanley-Reisner ring}.  As a Stanley-Reisner ring is associated with 
a simplicial complex, our $R/I$ is associated with a polyhedral complex 
$\Delta$ contained in the polyhedral cone $\R_{\geq 0}C \subset \R^d$. 
So we denote it by $K[\Delta]$. Our $K[\Delta]$ is a special case of 
the rings Stanley constructed from more general 
polyhedral complexes in \cite[\S4]{St2}, 
but still an interesting class. 
On the other hand, the sequentially Cohen-Macaulay property  
has become important in the theory of Stanley-Reisner rings, 
since it is closely related to 
{\it non-pure shellability} and {\it shifting} of simplicial complexes 
(c.f. \cite{St, D}).  Among other things,  
we show that the sequentially 
Cohen-Macaulay property of $K[\Delta]$ is a topological property of the 
``geometric realization" $|\Delta| \subset B$ of the complex $\Delta$. 
Here $B$ is a $(d-1)$-dimensional polytope which is the intersection 
of the cone $\R_{\geq 0}C$ and a hyperplane $H \subset \R^d$.

With the above notation, 
a squarefree module $M$ gives a $K$-constructible sheaf $M^+$ 
on $B$. For example, $K[\Delta]^+$ is the $K$-constant sheaf on the 
geometric realization $|\Delta| \subset B$. 
As shown in \cite{Y6}, we can connect the local cohomology of $M$ 
and the sheaf cohomology of $M^+$. More precisely, we have 
$[H_\m^{i+1}(M)]_0 \cong H^i(B;M^+)$ 
for all $i \geq 1$ and an exact sequence 
$0 \to [H_\m^0(M)]_0 \to M_0 \to H^0( B; M^+) \to [H_\m^1(M)]_0 \to 0$. 
We say $M \in \Sq$ is {\it regular}, if 
$[H_\m^0(M)]_0=[H_\m^1(M)]_0 = 0$. For all $M \in \Sq$, there 
is a unique $\overline{M} \in \Sq$ which is regular and  
$\overline{M}^+ \cong M^+$.  In \S5, we show that if $M$ is regular 
then the (sequentially) Cohen-Macaulay property of $M$ depends only on 
the sheaf $M^+$. This result can imply the main result of \S4. 

Assume that $R$ is normal and $I \subset R$ is a 
(not necessarily monomial) ideal. In \S6, generalizing a result of 
Herzog, Takayama and Terai~\cite{HTT},  
we show that if $R/I$ is Cohen-Macaulay and 
$\sqrt{I}$ is a monomial ideal, then $R/\sqrt{I}$ is Cohen-Macaulay again. 
(If $\sqrt{I}$ is not a monomial ideal, this statement does not hold. 
Otherwise, if $I \subset S:=K[x_1, \ldots, x_d]$ 
defines a set theoretic complete intersection subscheme of 
${\mathbb A}^d$, then $S/\sqrt{I}$ must be Cohen-Macaulay. 
This is clearly strange.) 

Let $S=K[x_1, \ldots, x_d]$ be a polynomial ring, and $I \subset S$ 
a monomial ideal. Recently, Takayama~\cite{T} showed that the range  
$\{ \, \ba \in \Z^d \mid [H_\m^i(S/I)]_\ba \ne 0 \, \}$  
is controlled by the degrees of minimal generators of $I$ 
(especially, when $H_\m^i(S/I)$ has finite length). 
In \S6, after giving a simple new proof of this result, we will 
(partially) extend it to affine semigroup rings. For example,  
we prove the following: Assume that $R$ is normal and simplicial. 
For $M \in \CgrR$, $H_\m^i(M)$ has finite length if and only if 
$H_\m^i(M)$ is $C$-graded.  
A similar result holds for arbitrary affine semigroup rings, 
but some modification is necessary.

\section{Irreducible resolutions over 
Cohen-Macaulay simplicial semigroup rings}
We use the same notation as in the introduction. 
Let $\R_{\geq 0} := \{ r \in \R \mid r \geq 0 \}$ be 
the set of non-negative real numbers, and  
$\bP:= \{ \, \sum \gamma_i \bc_i \mid \gamma_i \in \R_{\geq 0}, 
\bc_i \in C  \,  \} \subset \R^d$  the polyhedral cone spanned by $C$.   
Let $\cL$ be the set of (non-empty) faces of $\bP$. 
Note that $\bP$ itself and $\{ 0 \}$ belong to $\cL$. 
For $F \in \cL$, $P_F$ denotes the monomial ideal 
$(\bx^\bc \mid \bc \in C \setminus F)$ of $R$. Then $P_F$ is a prime ideal. 
Conversely, any monomial prime ideal of $R$ is of the form  $P_F$ for 
some $F \in \cL$. 

It is well-known that $\GrR$ has enough injectives 
(c.f. \cite[Chapter 11]{MS}). We denote the 
injective hull of $R/P_F$ in $\GrR$ by $\E(R/P_F)$. 
Assume that $M \in \GrR$ is indecomposable. Then $M$ is injective in $\GrR$ 
if and only if there are some $\ba \in \Z^d$ and 
$F \in \cL$ such that $M \cong \E(R/P_F)(-\ba)$. 
Recall that $M \in \GrR$ has a minimal injective resolution in $\GrR$, 
which is unique up to isomorphism. 
For a monomial prime ideal $P_F$, set 
$\mu_i(P_F, M):= \dim_{K(F)} (\Ext_R^i(R/P_F, M) \otimes R_{P_F})$, 
where $K(F)$ is the quotient field of $R/P_F$. 
It is well-known that if $\J$ is a minimal injective resolution of $M$ in 
$\GrR$, then we have $J^i \cong \bigoplus \E(R/P_F)^{\mu_i(P_F, M)}$ 
for all $i$ {\it as underlying $R$-modules} 
(i.e., if we forget the grading of the modules).  

We say an ideal $I \subset R$ is {\it irreducible} if every expression 
$I = I_1 \cap I_2$ of $I$ as an intersection of two ideals satisfies  
that $I = I_1$ or $I = I_2$. An irreducible ideal is always a primary 
ideal, while the converse is not true. In particular, the radical 
$\sqrt{I}$ of an irreducible ideal $I$ is a prime ideal. 
According to Miller~\cite{Mil}, we use the letter $W$ to denote 
a {\it monomial} irreducible ideal. 
Then we have $\sqrt{W} = P_F$ for some $F \in \cL$. 
In this case, $\dim (R/W)$ equals the dimension of $F$ as a polyhedral cone. 
When $R$ is a polynomial ring $K[x_1, \ldots, x_d]$, 
$I$ is a monomial irreducible ideal if 
and only if it is a complete intersection ideal 
generated by powers of variables.

For $M \in \GrR$, we say the submodule $\bigoplus_{\bc \in C} M_\bc$ 
is the {\it $C$-graded part} of $M$, and we denote it by $M_C$.  

\begin{prop}[{Miller \cite[\S2]{Mil}}]\label{Mil1}
We have the following.
\begin{itemize}
\item[(1)] Let $I \subset R$ be a monomial ideal. Then $I$ is 
irreducible  if and only if there are some $\ba \in C$ and 
$F \in \cL$ such that $R/I \cong [\E(R/P_F)(-\ba)]_C$. 
\item[(2)] The category $\CgrR$ has enough injectives. An indecomposable 
module $M$ is injective in $\CgrR$ if and only if $M \cong R/W$ for 
some monomial irreducible ideal $W \subset R$. 
\end{itemize}
\end{prop}

We say an injective resolution in $\CgrR$ is an 
{\it irreducible resolution}.

\begin{thm}[{Miller \cite[Theorem 2.4]{Mil}}]\label{Mil2}
Let $M \in \CgrR$ be a $C$-graded module, 
and $\J$ a minimal injective resolution of 
$M$ in $\GrR$. Then the $C$-graded part $[\J]_C$
of $\J$ is an irreducible resolution of $M$, and has finite length. 
\end{thm}

We say the irreducible resolution given in Theorem~\ref{Mil2} 
is a {\it minimal irreducible resolution}. 
Every $M \in \CgrR$ has a minimal irreducible resolution, 
and this is unique up to isomorphism. Any irreducible resolution is a direct 
sum of a minimal one and an exact sequence.  It is noteworthy that 
Helm and Miller~\cite{HM} gave an algorithm to compute minimal irreducible 
resolutions. 

Let $\I$ be a minimal irreducible resolution of $M \in \CgrR$. 
Then, for a monomial irreducible ideal $W$ and an integer $i \geq 0$, 
we have $\nu_i(W,M) \in \NN$ satisfying  
$$I^i \cong \bigoplus (R/W)^{\nu_i(W,M)}.$$ 
Note that 
\begin{equation}\label{mu & nu}
\mu_i(P_F,M) \geq \sum_{\sqrt{W}=P_F}\nu_i(W,M)
\end{equation}
for all $F \in \cL$. Set 
$\irrdim M := \max\{\, i \mid I^i \ne 0 \,\}$. Theorem~\ref{Mil2} 
states that $\irrdim M < \infty$ for all $M$. 
But, if $R$ is not a polynomial ring, then 
$\sup \{ \, \irrdim M \mid M \in \CgrR \, \} = \infty$. 
In fact, for a given integer $n$, we have 
$\irrdim (  \, (R/\m)(-\bc)  \, ) > n$  for sufficiently ``large" $\bc \in C$.  

\begin{lem}\label{associated prime}
Let $M \in \CgrR$. If $\nu_0(W,M) \ne 0$, then $\sqrt{W}$ is an associated 
prime of $M$. Conversely, if $P_F$ is an associated prime of $M$, 
then there is a monomial irreducible ideal $W$ with 
$\sqrt{W} = P_F$ and $\nu_0(W,M) \ne 0$. 
In particular, $\dim M = \max \{ \, \dim (R/W) \mid \nu_0(W,M) \ne 0 \, \}$. 
\end{lem}

\begin{proof} 
The first assertion follows from \eqref{mu & nu}. 
Since $M$ is a submodule of $I^0 := \bigoplus (R/W')^{v_0(W',M)}$ and 
an irreducible ideal is a primary ideal, we have 
$\Ass(M) \subset \Ass(I^0) = \{ \sqrt{W'} \mid v_0(W',M) \ne 0 \}$. 
So the second statement follows. The last assertion is easy.
\end{proof}

For further information on irreducible resolutions, consult 
\cite[Chapter~11]{MS}. 

\bigskip

We say $R$ is {\it simplicial} if there are 
$\be_1, \cdots, \be_d \in C$ such that $\bP = 
\{ \, \sum_{i=1}^d \gamma_i \be_i \mid \gamma_i \in \R_{\geq 0} \, \}$. 
In this case, we have $C=\Z^d \cap \{ \, \sum_{i=1}^d 
\alpha_i \be_i \mid \alpha_i \in \Q_{\geq 0} \,\}$, and 
$\{\be_1, \cdots, \be_d \}$ is a basis of $\R^d$ (and $\Q^d$). 
A polynomial ring $K[x_1, \ldots, x_d]$ is a primary example of 
a simplicial semigroup ring. 

\begin{lem}\label{R/W is CM}
Assume that $R$ is Cohen-Macaulay and 
simplicial. Then $R/W$ is Cohen-Macaulay for all monomial 
irreducible ideal $W$.  
\end{lem}

\begin{proof}
By Proposition~\ref{Mil1} (1), 
there is some $\ba =\sum_{i=1}^d \alpha_i \be_i \in C$ 
($\alpha_i \in \Q_{\geq 0}$) such that $[\E(R/P_F)(-\ba)]_C \cong R/W,$ 
where $F$ is the face of $\bP$ with $P_F = \sqrt{W}$. 
We may assume that $F$ is spanned by $\be_1, \be_2, \ldots, \be_n$ 
as a polyhedral cone. 
Set $H:= C \cap F$ to be a submonoid of $C$, and set 
$$U := \{ \, \sum_{i=1}^d \beta_i \be_i \in C \mid 
\text{$\beta_i \geq 0$ for all $i$ and 
$\beta_i \leq \alpha_i$ for all $i > n$}  \, \}.$$
Then $U$ is an $H$-{\it ideal} in the sense of \cite{GSW}, 
that is, $H + U \subset U$. 
(In \cite{GSW}, they assumed that $U$ is contained in the group generated 
by $H$. But their results hold without this assumption.) 
Note that $A:= K[H] = \bigoplus_{\bc \in H} K \bx^\bc$  is a simplicial 
affine semigroup ring. Clearly, $A \cong R/P_F$. But here we regard 
$A$ as a subring of $R$ via the inclusion $H \subset C$. 
As an $A$ module, $R/W$ is isomorphic to $K[U]$, and it is finitely generated. 
Since $R=K[C]$ is Cohen-Macaulay, $U$ satisfies the condition (iii) 
of \cite[Theorem~2,2]{GSW} as an $H$-ideal. Thus $R/W$ is a 
Cohen-Macaulay module over $A$ (thus over $R$). 
\end{proof}

\begin{exmp}\label{square}
(1) In the non-simplicial case, a monomial irreducible
ideal need not be Cohen-Macaulay, even if the ring is normal. 
Consider the affine semigroup 
$C$ in $\Z^3$ generated by $(0,0,1)$, $(0, 1, 1)$, 
$(1,1,1)$ and $(1,0,1)$. Then the semigroup ring $R$ of $C$ is 
normal but not simplicial. We denote the monomials in $R$ with 
degrees $(0,0,1)$, $(0, 1, 1)$, $(1,1,1)$ and $(1,0,1)$ by 
$x$, $y$, $z$ and $w$ respectively. Note that $R = K[x,y,z,w]/(xz-yw)$.  
Set $W = (y^2, yz, z^2)$.  Then $W$ is an irreducible ideal,  
since $R/W \cong [\E(R/(y,z))(-(1,1,1))]_C$. 
But computation by {\it Macaulay 2} shows that $R/W$ is not Cohen-Macaulay.  

(2) If $R$ is normal (but not necessarily simplicial), 
then $R/P_F$ is a normal semigroup ring (in particular, Cohen-Macaulay) 
for any monomial prime ideal $P_F$. 
If $R$ is simplicial and Cohen-Macaulay, then so is $R/P_F$ 
by Lemma~\ref{R/W is CM}. On the other hand, if $R$ is Cohen-Macaulay 
but not simplicial, then $R/P_F$ need not be Cohen-Macaulay. 
The following example is due to Professor Ngo Viet Trung.  
The affine semigroup ring 
$R:= k[s^4, s^3t, st^3, t^4, s^4u, t^4u]$ is Cohen-Macaulay. 
But $P=(s^4u, t^4u)$ is a prime ideal and $R/P \cong k[s^4, s^3t, st^3, t^4]$ 
is not Cohen-Macaulay. 
\end{exmp}

When $R$ is Cohen-Macaulay, $\can$ denotes the canonical module of $R$. 

\begin{thm}\label{depth}
If $R$ is Cohen-Macaulay and simplicial, and  $M \in \CgrR$, then we have 
\begin{equation}\label{depth-eq}
\depth_R M = \min \{ \, \dim (R/W) + i \mid  \nu_i(W,M) \ne 0 \, \}.
\end{equation}
\end{thm}

\begin{proof}
Let $\I: 0 \to I^0 \stackrel{d^0}{\to} I^1  \stackrel{d^1}{\to} \cdots$ 
be a minimal irreducible resolution of $M$, and set 
$\Omega^i(M) := \ker(d^i)$. Of course, $\Omega^0(M) = M$. 
Set $\delta$ to be the right hand side of \eqref{depth-eq}. 

To prove $\depth_R M \geq \delta$, we will show that 
$\depth (\Omega^i(M)) \geq \delta -i$ for all $i$ by 
backward induction on $i$.  
(Here we set the depth of the 0 module to be $+ \infty$.) 
If $i \geq d$, there is nothing to prove. 
Assume that $\depth_R (\Omega^{i+1}(M)) \geq  \delta -i-1$. 
Consider the short exact sequence 
\begin{equation}\label{short exact}
0 \to \Omega^i(M) \to I^i \to \Omega^{i+1}(M) \to 0.
\end{equation}
Recall that $I^i = \bigoplus (R/W)^{\nu_i(W,M)}$ and each $R/W$ is 
Cohen-Macaulay by Lemma~\ref{R/W is CM}.  
By the assumption, if  $\nu_i(W,M) \ne 0$ then $\dim (R/W) \geq \delta - i$.  
So $\depth_R I^i \geq \delta - i$. By \eqref{short exact}, 
we have $\depth_R (\Omega^i(M)) \geq \delta - i$. 

We can take $n \in \NN$  such that  $\nu_n(W,M) \ne 0$ 
for some $W$ with $\dim (R/W) = \delta-n$. 
To prove $\depth_R M = \delta$, we will show that 
$\depth_R (\Omega^i(M)) = \delta -i$ for all $i \leq n$ by 
backward induction on $i$. 
Since $\sqrt{W}$ is an associated prime of $\Omega^n(M)$ 
and $\dim (R/\sqrt{W}) = \delta - n$, we have 
$\depth_R (\Omega^n(M)) \leq \delta - n$. (In fact, we have 
$H_\m^{\delta - n}(\Omega^n(M))^\vee \cong 
\Ext_R^{d-\delta +n} (\Omega^n(M),\can) \ne 0$ by an argument similar to 
\cite[Theorem~8.1.1]{BH}.) 
But we have seen that $\depth_R (\Omega^n(M)) \geq \delta -n$. 
So $\depth_R (\Omega^n(M)) = \delta  - n$. Assume that 
$\depth_R (\Omega^{i+1}(M)) =  \delta -i-1$. 
Since $\depth I^i \geq \delta - i$, we have 
$\depth_R (\Omega^i(M)) =  \delta -i$ by \eqref{short exact}. 
\end{proof}

The following is an easy consequence of the theorem. 

\begin{cor}\label{CM}
Assume that $R$ is Cohen-Macaulay and simplicial.  
Then $M \in \CgrR$ is a Cohen-Macaulay module of dimension $p$ if 
and only if $\nu_i(W,M) \ne 0$ implies $p-i \le \dim(R/W) \leq p$ 
for all $i$. 
\end{cor}

We can also characterize Serre's condition $(S_n)$ in our context. 

\begin{thm}\label{S_n}
Assume that $R$ is Cohen-Macaulay and simplicial. 
If $I \subset R$ is a monomial 
ideal with $\dim (R/I) = p$, then  
the following are equivalent for an integer $n \geq 2$. 
\begin{itemize}
\item[(a)] $R/I$ satisfies Serre's condition $(S_n)$. 
\item[(b)] If $i \leq n-1$ and $\dim (R/W) < p-i$, then 
$\nu_i(W, R/I) =0$.  
\end{itemize}
\end{thm}

To prove the theorem, we need the following 
(more or less) well-known facts. 
Note that these facts hold in much wider context 
(e.g., over a noetherian local ring admitting a dualizing complex).

\begin{lem}\label{Ext^j}
For a monomial ideal $I \subset R$ and a module $M \in \grR$,  
we have the following. \begin{itemize}
\item[(1)] We always have $\dim (\Ext_R^j(M,\can)) \leq d-j$. 
And the equality holds if and only if there is an associated prime $P$ of 
$M$ with $\dim (R/P) = d-j$. 
\item[(2)] If $R/I$ satisfies the  $(S_n)$ condition for some $n \geq 2$, 
then all associated primes of $R/I$ have the same dimension. 
\item[(3)] Let $n \geq 2$. Then $R/I$ satisfies $(S_n)$ 
if and only if 
$\dim (\Ext^j_R(R/I, \can)) \leq d-j-n$ for all $j > d- \dim (R/I)$. 
\end{itemize}
\end{lem}

\begin{proof}
(1) Essentially same as \cite[Theorem~8.1.1]{BH}.

(2) See \cite[Remark~2.4.1]{Hart}. 

(3) Follows from (2) and the local duality. 
\end{proof}

The parts (2) and (3) of Lemma~\ref{Ext^j} do not hold for a module $M$. 
So Theorem~\ref{S_n} only concerns ideals. But if all minimal primes of 
$M \in \CgrR$ have the same dimension, then Theorem~\ref{S_n} holds for $M$. 

\bigskip

\noindent{\it Proof of Theorem~\ref{S_n}.} Set $A := R/I$. 
As before, let $\I: 0 \to I^0 \stackrel{d^0}{\to} I^1  \stackrel{d^1}{\to} 
\cdots$ be a minimal irreducible resolution of $A$, and set 
$\Omega^i(A) := \ker(d^i)$.

Assume that $A$ satisfies the condition (b). 
Then $\Ext_R^j(I^i, \can) = 0$ for all $i < n$ and $j > d-p+i$ by 
Lemma~\ref{R/W is CM}. The short exact sequence \eqref{short exact} yields 
$$\Ext_R^j(\Omega^i(A), \can) \cong \Ext_R^{j+1}(\Omega^{i+1}(A), \can)$$
for all $i < n-1$ and $j > d-p+i$. Using this, we get 
$$\Ext_R^j(A, \can) \cong \Ext_R^{j+n-1}(\Omega^{n-1}(A), \can)$$
for all $j > d-p$. In this case, the dimension of 
any associated prime of $\Omega^{n-1}(A)$ is at least $p-n+1$ by (b), 
and we have $\dim (\Ext_R^{j+n-1}(\Omega^{n-1}(A), \can)) \leq d-j-n$ 
by Lemma~\ref{Ext^j} (1) (note that $j+n-1 > d-p+n-1$). 
Thus $A$ satisfies the $(S_n)$ condition by 
Lemma~\ref{Ext^j} (3).
 
Conversely, assume that $A$ satisfies $(S_n)$ for some $n \geq 2$ but  
$\nu_i(W,A) \ne 0$ for some $i \leq n-1$ and some $W$ with $\dim (R/W) \leq 
p-i-1$. Set $P = \sqrt{W}$. Then $\dim A_P \geq p -(p-i-1) = i+1$. 
But, we have $\Ext_R^i(R/P, A) \ne 0$ by \eqref{mu & nu}. 
Hence $\depth_{R_P} A_P \leq i < \min\{ n, \dim A_P \}$. 
It contradicts $(S_n)$.   
\qed

\begin{dfn}\label{gCM}
We say $M \in \grR$ is a {\it  generalized Cohen-Macaulay module } 
if $H_\m^i(M)$ has finite length 
(i.e., $\dim_K H_\m^i(M)<\infty$) for all $i <  \dim M$. 
We can also define this concept 
over a noetherian local ring in a similar way. 
\end{dfn}

By Lemma~\ref{Ext^j} (1), all minimal primes of a 
generalized Cohen-Macaulay module have the same dimension. 

\begin{prop}\label{gCM by resol}
Assume that $R$ is Cohen-Macaulay and simplicial. 
For $M \in \CgrR$ with $\dim M = p$,  
$M$ is a generalized Cohen-Macaulay module if and only if 
$\nu_i(W,M) = 0$ for all $i$ and all $W$ 
with $0 < \dim (R/W) < p-i$. 
\end{prop}

\begin{proof}
(Sufficiency) By argument similar to the proof of Theorem~\ref{depth}, 
$H_\m^j(\Omega^i(M))$ has finite length for all $j < p-i$. 
Since $\Omega^0(M) = M$, we are done. 

(Necessity)  We will prove the contrapositive. 
Assume that $\nu_n(W,M) \ne 0$ for some $n \geq 0$ and some 
$W$ with $0 < \dim (R/W) < p-n$. 
We may assume that $n$ is minimum among such integers, and 
set $q := \dim (R/W)$. We can prove that 
$H_\m^{q + n -i}(\Omega^i(M))$ does not have finite length  for all 
$i \leq n$ by backward induction. 
Since $\Omega^0(M) = M$ and $q+n < p$, $M$ is not generalized Cohen-Macaulay. 
\end{proof}

The notion of {\it Buchsbaum modules} is an intermediate concept between 
Cohen-Macaulay modules and generalized Cohen-Macaulay modules.
\cite{SV} is a nice introduction of this theory. 

\begin{prop}\label{Bbm}
Assume that $R$ is Cohen-Macaulay and simplicial. 
Let $M \in \CgrR$ be a $C$-graded module of dimension $p$. If  
$\nu_i(W,M) = 0$ for all $i \geq 0$ and all monomial irreducible ideal $W$ 
with $W \ne \m$ and $\dim (R/W) < p-i$, then $M$ is Buchsbaum. 
\end{prop}

\begin{proof}
By an argument similar to the above, 
we can see that $[H_\m^i(M)]_\ba = 0$ for all $i < p$ and all 
$0 \ne \ba \in \Z^d$. So the assertion follows from \cite[Corollary~4.6]{Y6}. 
\end{proof}

\section{Squarefree modules over a normal semigroup ring}
For the results in the previous section, the assumption that 
$R$ is simplicial is really necessary. But when we consider radical monomial 
ideals in a normal semigroup ring, we can remove this assumption. 

{\it Throughout this section, we assume that $R$ is normal.} 

For a point $u \in \bP \, (=\R_{\geq 0} C)$,  we always have a unique face 
$F \in \cL$ whose relative interior contains $u$. 
In this notation, we denote $s(u) = F$.

\begin{dfn}[\cite{Y2}]\label{sq}
Assume that $R$ is normal. We say an $R$-module $M$ 
is {\it squarefree}, if $M \in \CgrR$ and the multiplication map 
$M_\ba \ni y \mapsto \bx^\bb y 
\in M_{\ba + \bb}$ is bijective for all 
$\ba, \bb \in C$ with $s(\ba+\bb) = s(\ba)$. 
\end{dfn} 

For a monomial ideal $I \subset R$, $I$ (or $R/I$) is a squarefree module 
if and only if $I$ is a radical ideal (i.e., $I=\sqrt{I}$).  
Since the canonical module $\can$ of $R$ is isomorphic to the ideal 
$( \, \bx^\bc \mid \text{$\bc \in C$ with $s (\bc) = \bP$} \, )$, 
it is also squarefree. It is easy to check that if $M$ is squarefree, 
we have $\dim_K M_\ba = \dim_K M_\bb$ for all $\ba, \bb \in C$ 
with $s(\ba) = s(\bb)$.

Let us recall basic properties of squarefree modules. 
See \cite{Y2} for detail. 
$\Sq$ denotes the full subcategory of $\GrR$ consisting of 
squarefree modules. Then $\Sq$ is closed under kernels, cokernels 
and extensions in $\GrR$. 
Hence $\Sq$ is an abelian category. Moreover, it admits enough 
projectives and injectives, and an indecomposable injective object 
is isomorphic to $R/P_F$ for some $F \in \cL$.  

If $M \in \grR$ and $N \in \GrR$, then the $R$-module $\Hom_R(M,N)$
has a natural $\Z^d$-grading with 
$\Hom_R(M,N)_\ba = \Hom_{\GrR}(M,N(\ba))$. 
Similarly, $\Ext_R^i(M,N)$ has a natural  $\Z^d$-grading.

\begin{lem}[\cite{Y2}]
Let $M$ be a squarefree $R$-module, and $\I$ its minimal irreducible 
resolution. Then we have the following.  
\begin{itemize}
\item[(1)] If $\nu_i(W,M) \ne 0$ for some $i \geq 0$, 
then $W$ is a prime ideal. Moreover, $\I$ is a minimal 
injective resolution of $M$ in $\Sq$. 
\item[(2)] $\dim I^i > \dim I^{i+1}$ for all $i$. 
In particular, $\irrdim M \leq d$. 
\item[(3)] $\Ext^i_R(M, \can)$ is squarefree for all $i$. 
\end{itemize}
\end{lem}

Let $M$ be a squarefree module. For each $F \in \cL$, 
take some $\bc(F) \in C \cap \relint(F)$ (i.e., $s(\bc(F)) =F$). 
If $F, G \in \cL$ and $G \supset F$, \cite[Theorem~3.3]{Y2} gives a 
$K$-linear map $\varphi^M_{G, F}: M_{\bc(F)} \to 
M_{\bc(G)}$.  These maps satisfy 
$\varphi^M_{F,F} = \Id$ and 
$\varphi^M_{H, G} \circ \varphi^M_{G, F} = \varphi^M_{H,F}$ for all 
$H \supset G \supset F$.

\begin{thm}[{\cite[Theorem~4.15]{Y2}}]\label{Ext_F} 
For $M\in \Sq$ and $F \in \cL$ with $\dim F = t$, we have 
\begin{equation}\label{Ext}
\nu_i(P_F,M) = \dim_K [\Ext^{d-i-t}_R (M, \can)]_{\bc(F)}.
\end{equation} 
\end{thm}

Since $\Ext^{d-i-t}_R (M, \can)$ is squarefree, 
the value of the right side of the equality \eqref{Ext} 
does not depend on the choice of $\bc(F)$.

\medskip 

Propositions~\ref{CMsq}, \ref{Snsq} and \ref{gCMsq}, which 
follow from Theorem~\ref{Ext_F},  
are the squarefree modules version of results in the previous section. 
While (the latter part of) Proposition~\ref{CMsq} has been stated in 
\cite{Y2}, we state it here for the reader's convenience. 
Since $R/P_F$ is always Cohen-Macaulay whenever $R$ is normal, 
we can also prove these results by arguments similar to 
the previous section. 

\begin{prop}[\cite{Y2}]\label{CMsq}
If $M \in \Sq$,  
then $\depth_R M = \min \{ \, \dim (R/P_F) + i \mid 
\nu_i(P_F,M) \ne 0 \, \}.$ In particular, 
$M$ is a Cohen-Macaulay module of dimension $p$ if 
and only if $\nu_i(P_F,M) = 0$ for all $i$ and all $P_F$ with 
$\dim (R/P_F) \ne p -i$ .
\end{prop}

\begin{prop}\label{Snsq}
If $I \subset R$ is a radical monomial ideal with $\dim (R/I) = p$, then  
the following are equivalent for an integer $n \geq 2$. 
\begin{itemize}
\item[(a)] $R/I$ satisfies Serre's condition $(S_n)$. 
\item[(b)] If $i \leq n-1$ and $\dim (R/P_F) \ne p-i$, then 
$\nu_i(P_F, R/I) =0$.  
\end{itemize}
\end{prop}

For squarefree modules, we can prove the converse of Proposition~\ref{Bbm} 
by virtue of \cite[Corollary~4.6]{Y6}.

\begin{prop}\label{gCMsq}
Assume that $R$ is normal. Let $M$ be a squarefree $R$-module of dimension 
$p$. Then $M$ is Buchsbaum, if and only if it is generalized Cohen-Macaulay,
if and only if $\nu_i(P_F, M)=0$ for all $i$ and all monomial prime ideal 
$P_F$ with $\dim (R/P_F) \ne 0, p-i$.   
\end{prop}

The {\it linear strand} of a minimal free resolution of a finitely 
generated graded module over a polynomial ring is an important 
notion introduced by Eisenbud (c.f. \cite{Ei}). 
Here we introduce the analog of this concept for an irreducible 
resolution of squarefree modules. (Miller also studied this concept 
implicitly. See \cite[\S3]{Mil}.)
 
Let $M$ be a squarefree $R$-module, and $\I$ its minimal irreducible 
resolution. For each $l \in \NN$, 
we define the {\it $l$-linear strand} 
$\lin_l(M)$ of (the minimal irreducible resolution of) $M$ as follows: 
The term $\lin_l(M)^i$ of cohomological degree $i$  
is $$\bigoplus_{\dim(R/P_F) = l-i} (R/P_F)^{\nu_i(P_F,M)},$$
which is a direct summand of $I^i = \bigoplus (R/P_F)^{\nu_i(P_F,M)}$, 
and the differential $\lin_l(M)^i \to \lin_l(M)^{i+1}$ 
is the corresponding component of the differential 
$I^i \to I^{i+1}$ of $\I$.  By Proposition~\ref{CMsq}, 
$M$ is a Cohen-Macaulay module of dimension $p$ 
if and only if $\I \cong \lin_p(M)$. 
Set $\lin(M) := \bigoplus_{0 \leq l \leq d} \lin_l(M)$. 

Using spectral sequence argument, we can construct $\lin (M)$ from  
a (not necessarily minimal) irreducible resolution $\J$ of $M$. 
For each $i \in  \NN$, 
let $J^i \cong \bigoplus_{j \geq 0} J^{i,j}$ be a decomposition 
with $\dim (J^{i,j}) = j$ or $J^{i,j} = 0$. Consider the filtration 
$\J = F_0\J \supset F_1 \J \supset \cdots \supset F_d \J =0$ with  
$F_p J^i = \bigoplus_{j \leq d-p} J^{i,j}$, and  
the associated spectral sequence $\{E_r^{*,*},d_r \}$. 
Since $F_pJ^i$ is a  direct summand of $J^i$, we have 
$E_0^{p,q} = (F_p J^{p+q} / F_{p+1} J^{p+q}) \cong J^{p+q, d-p}$ and 
$J_0^t := \bigoplus_{p+q=t}E_0^{p,q} \cong J^t$.  
The maps $d_0^{p,q} : E_0^{p,q} \to E_0^{p,q+1}$ make $\J_0$ 
a cochain complex.  On the other hand, 
we always have a decomposition $\J = \I \oplus \T$ such that $\I$ is minimal 
and $\T$ is exact. If we identify $J_0^t$ with 
$J^t = I^t \oplus T^t$, the differential $d_0$ of $\J_0$ is given by 
$(0, d_{\T})$. As before, let $I^i \cong \bigoplus_{j \geq 0} I^{i,j}$ 
be a decomposition with $\dim (I^{i,j}) = j$ or $I^{i,j} = 0$. 
Then we have $E_1^{p,q}  \cong F_p I^{p+q} /F_{p+1} I^{p+q} \cong 
I^{p+q, d-p}$ and $J_1^t := \bigoplus_{p+q=t}E_1^{p,q} \cong I^t$. 
The maps $d_1^{p,q}: E_1^{p,q} \to E_1^{p+1,q}$ make $\J_1$ 
a cochain complex, which is isomorphic to $\lin(M)$.  

A minimal injective resolution of $\can$ in $\GrR$ 
is of the form \begin{equation}\label{dual}
\I: 0 \to I^0 \to I^1  \to \cdots \to I^d \to 0,
\end{equation}
$$I^i = \bigoplus_{\substack{F \in \cL \\ \dim F = d-i}} \E(R/P_F)$$
and the differential 
is composed of the maps $(\pm1) 
\operatorname{nat} : \E(R/P_F) \to \E(R/P_{F'})$ for $F,F' \in \cL$ 
with $F \supset F'$ and $\dim F = \dim F'+1$, 
where $\operatorname{nat} : \E(R/P_F) \to \E(R/P_{F'})$ is 
induced by the natural surjection $R/P_F \to R/P_{F'}$, and 
the sign $\pm$ is given by an {\it incidence function} 
(c.f. \cite[\S6.2]{BH}) of $\bP$. 

Let $D^b(\Sq)$ be the bounded derived category of $\Sq$. 
By \cite[Lemma~2.3]{Y6}, we have $D^b(\Sq) \cong D^b_{\Sq}(\GrR)$, 
so we will identify these categories. 
For a complex $\M$ and an integer $n$, let $\M[n]$ be the $n^{\rm th}$ 
translation of $\M$. That is, $\M[n]$ is a complex with $M^i[n] = M^{i+n}$. 
Since $\Ext^i_R(M,\can) \in \Sq$ for all $M \in \Sq$, 
$$\Du(-):= \RHom_R(-, \can)$$ gives a contravariant functor 
from $D^b(\Sq)$ to itself satisfying $\Du \circ \Du \cong \Id_{D^b(\Sq)}$. 

For $M \in \GrR$, we say the submodule $\bigoplus_{\bc \in C} M_\bc$ 
is the {\it $C$-graded part} of $M$. 
Since $R$ is normal now, the $C$-graded part of $\E(R/P_F)$ is 
isomorphic to $R/P_F$ (c.f. \cite[Chapter~13]{MS}). 
Moreover, we have the following.  

\begin{lem}\label{C-part} 
Assume that $R$ is normal. For $M \in \Sq$ and $F \in \cL$,  
the $C$-graded part of $\Hom_R(M, \E(R/P_F))$ 
is isomorphic to $(M_{\bc(F)})^* \otimes_K (R/P_F)$. 
\end{lem}

\begin{proof}
When $R$ is a polynomial ring, this result was given in 
\cite[Lemma~3.20]{Y1}. The general case can be proved by essentially 
same argument. But we give a precise proof here for the reader's 
convenience.

Let $M'$ be the submodule of $M$ generated by $M_{\bc(F)}$. 
By the injectivity of $\E(R/P_F)$, 
$f \in \Hom_{\GrR} (M', \E(R/P_F)(\bc))$ for $\bc \in C$ 
can be extended to $M \to \E(R/P_F)(\bc)$. On the other hand, 
if $\bc' \in C \setminus F$, then $[\E(R/P_F)]_{\bc+\bc'}=0$. 
Hence, for $0 \ne g \in \Hom_{\GrR}(M, \E(R/P_F)(\bc))$, 
there is some $\bc' \in F$ such that $g(y) \ne 0$ for some $y \in M_{\bc'}$. 
In this situation, we have $g(\bx^{\bc(F)}y)= \bx^{\bc(F)} g(y) \ne 0$. 
Since $M$ is squarefree and $s(\bc(F))=s(\bc(F) + \bc')$, 
there is $z \in M_{\bc(F)}$ such that $\bx^{\bc'} z = \bx^{\bc(F)} y$. 
Clearly, $g(z) \ne 0$. Hence the restriction $g |_{M'}$ of 
$g \in \Hom_{\GrR}(M, \E(R/P_F)(\bc))$ to $M'$ is 0 if and only if $g=0$. 
Combining these observations, we have $$[\Hom_R (M', \E(R/P_F))]_C 
\cong [\Hom_R(M, \E(R/P_F))]_C.$$ 

Since $\ann(y) \subset P_F$ for all $0 \ne y \in M_{\bc(F)}$, 
$f \in \Hom_K(M_{\bc(F)}, \E(R/P_F)(\bc))$ can be extended to 
$\tilde{f} \in \Hom_{\GrR}(M', \E(R/P_F)(\bc))$. 
Note that $\E(R/P_F)_{\bc + \bc(F)} \ne 0$ if and only if 
$\E(R/P_F)_{\bc + \bc(F)} \cong K$ if and only if $\bc \in F$. 
Hence we have $[\Hom_R (M', \E(R/P_F))]_C \cong 
(M_{\bc(F)})^* \otimes_K R/P_F$ as graded $K$-vector spaces. 
But it is easy to see that this is also an isomorphism of $R$-modules. 
\end{proof}

\begin{lem}\label{explicit form}
Assume that $R$ is normal. If $M \in \Sq$, then $\Du(M)$ is 
quasi-isomorphic to the complex 
$D^{\bullet} : 0 \to D^0 \to D^1 \to \cdots \to D^d \to 0$ with 
\begin{equation}\label{D(M)}
D^i = \bigoplus_{\substack{F \in \cL \\ \dim F = d-i}} 
(M_{\bc(F)})^* \otimes_K (R/P_F).
\end{equation}
Here the differential is the sum of the maps 
$$(\pm\varphi_{F,F'}^M)^* \otimes \nat
: (M_{\bc(F)})^* \otimes_K R/P_F \to (M_{\bc(F')})^* \otimes_K 
R/P_{F'}$$
for  $F,F' \in \cL$ with $F \supset F'$ and $\dim F = \dim F'+1$, 
and $\nat$ denotes the natural surjection $R/P_F \to R/P_{F'}$. 

Moreover, if $\M \in D^b(\Sq)$, then $\Du(\M)$ is quasi-isomorphic 
to the total complex of the double complex 
$W^{i,j}(\M) = D^i(M^{-j})$.  
Here the differential $W^{i,j}= D^i(M^{-j}) \to W^{i+1,j} = 
D^{i+1}(M^{-j})$ is the $i^{\rm th}$ differential of 
$D^\bullet(M^{-j})$ and the differential 
$$W^{i,j}=  \bigoplus_{\substack{F \in \cL \\ \dim F = d-i}} 
(M^{-j}_{\bc(F)})^* \otimes_K (R/P_F) \longrightarrow  W^{i,j+1} = 
\bigoplus_{\substack{F \in \cL \\ \dim F = d-i}} 
(M^{-j-1}_{\bc(F)})^* \otimes_K (R/P_F)$$ is induced by 
the differential $\partial_{\M}^{-j-1} : M^{-j-1} \to M^{-j}$ of $\M$. 
\end{lem}

\begin{proof}
When $R$ is a polynomial ring, the assertion was given in 
\cite[\S3]{Y5}. The general case can be proved by a similar argument 
to \cite{Y5}. Here we only remark that  
$\RHom_R(\M, \can)$ is quasi-isomorphic to $\Hom_R^\bullet(\M,\I)$ 
where $\I$ is the minimal injective resolution of $\can$ 
described in \eqref{dual}, 
and the isomorphism \eqref{D(M)} clearly follows from Lemma~\ref{C-part}. 
\end{proof}

\noindent{\bf Convention.} In the sequel, as an explicit complex, 
$\Du(\M)$ for $\M \in D^b(\Sq)$ means the complex described in 
Lemma~\ref{explicit form}. 

\medskip

The next result refines Theorem~\ref{Ext_F}.

\begin{thm}\label{lin(M)}
Assume that $R$ is normal. If $M \in \Sq$, then we have $\lin_i(M) \cong 
\Du(\Ext_R^{d-i}(M,\can))[d-i]$ for all $i$.
\end{thm}

\begin{proof}
Note that $\Du\circ\Du(M)$ is a complex of injective objects in $\Sq$, and 
quasi-isomorphic to $M$. Hence it is a (non-minimal) irreducible resolution 
of $M$. Set $\J := \Du \circ \Du(M)$ and $\N := \Du(M)$. 
Recall that $\J=\Du(\N)$ is the total complex of the double complex 
$W^{\bullet \bullet}$ defined in Lemma~\ref{explicit form}. 
We use the same notation as in the spectral sequence construction of $\lin(M)$. Recall that $J^i \cong J_0^i$ for all $i$. 
The construction of $\J_0$ cancels the horizontal differential 
of $W^{\bullet \bullet}$ (i.e., the differential which comes from 
that of $\Du(N^i)$). Hence the differential of $\J_0$ is just induced 
by that of $\N$. So if we set $\bH(\N)$ to be the complex 
such that $\bH(\N)^i = H^i(N)$ for all $i$ 
and the differential maps are zero, 
then $\J_1$ is isomorphic to $\Du(\bH(\N))$. 
Since the $i^{\rm th}$ cohomology of $\bH(\N)$ is 
$H^i(N) = \Ext_R^i(M, \can)$, the assertion follows. 
\end{proof}

\begin{dfn}[{Stanley, \cite{St}}] Let $M \in \grR$. 
We say $M$ is {\it sequentially Cohen-Macaulay} if there is a 
finite filtration 
$$0 = M_0 \subset M_1 \subset \cdots \subset M_r = M$$ 
of $M$ by graded submodules $M_i$ satisfying the following conditions.  
\begin{itemize}
\item[(a)] Each quotient $M_i/M_{i-1}$ is Cohen-Macaulay.
\item[(b)] $\dim(M_i/M_{i-1}) < \dim (M_{i+1}/M_i)$ for all $i$. 
\end{itemize}
\end{dfn}

\begin{cor}\label{seqCM}
Let $M \in \Sq$. 
Then $M$ is sequentially Cohen-Macaulay if and only if 
$\lin(M)$ is acyclic (i.e., $H^i(\lin(M))=0$ for all $i \ne 0$). 
\end{cor}

\begin{proof}
For $0 \ne N \in \Sq$, $H^i(\Du(N)) = 0$ for all $i \ne n$ 
if and only if $N$ is Cohen-Macaulay and $\dim N = d-n$. 
By Theorem~\ref{lin(M)}, $\lin(M)$ is acyclic if and only if 
$\Ext_R^{d-i}(M,\can)$ is a Cohen-Macaulay module of dimension $i$ 
for all $i$. So the assertion follows from \cite[Theorem~III.2.11]{St}. 
\end{proof}

\begin{rem}
Assume that $R$ is a polynomial ring. It is easy to see that $\lin(M)$ of 
$M \in \Sq$ is acyclic if and only if the Alexander dual of $M$ is 
{\it componentwise linear} (see, for example, \cite{R1, Y5} for this concept). 
So Corollary~\ref{seqCM} generalizes \cite[Theorem~4.5]{R1}. 
\end{rem}

\section{Sequentially Cohen-Macaulay face rings}
In this section, we always assume that $R$ is normal. 

We have a hyperplane $H \subset \R^d$ such that $B:= H \cap \bP$ 
is a $(d-1)$-dimensional polytope. 
Clearly, $B$ has the topology induced by the Euclidean space 
$\R^d$, and $B$ is homeomorphic to a closed ball of dimension $d-1$. 
For a face $F \in \cL$, set $|F|$ to be the relative interior of $F \cap H$.  
Regarding $\cL$ as a partially ordered set by inclusions, 
we say $\Delta \subset \cL$ is an {\it order ideal}, if $\Delta$ is 
a {\it non-empty} subset such that; 
$F \in \Delta$, $G \in \cL$ and $F \supset G$ $\Rightarrow$ $G \in \Delta$. 
If $\Delta$ is an order ideal,  
then $|\Delta|  := \bigcup_{F \in \Delta} |F|$ is a closed subset of $B$, and 
$\bigcup_{F \in \Delta} |F|$ is a {\it regular cell decomposition}
(c.f. \cite[\S 6.2]{BH}) of $|\Delta|$.  
Up to homeomorphism, (the regular cell decomposition of) $|\Delta|$ does 
not depend on the particular choice of the hyperplane $H$. 
The dimension  $\dim |\Delta|$ of $|\Delta|$ is given by 
$ \max \{ \, \dim |F| \mid F \in \Delta \, \}$. Here $\dim |F|$ denotes 
the dimension of $|F|$ as a cell (we set $\dim \emptyset = -1$), 
that is, $\dim |F| = \dim F -1$. 

We assign an order ideal $\Delta \subset \cL$ to the  
ideal $I_\Delta := (\, \bx^\bc \mid \text{$\bc \in C$ and 
$s(\bc) \not \in \Delta$} \, )$ of $R$. (For the definition of $s(\bc)$, 
see the beginning of \S3.) Note that $I_\Delta$ is a radical ideal, and 
any radical monomial ideal of $R$ is of the form $I_\Delta$ for some 
$\Delta$. Set $K[\Delta] := R/I_\Delta$. Clearly, 
$$
K[\Delta]_\bc \cong 
\begin{cases}
K & \text{if $\bc \in C$ and $s(\bc) \in \Delta$,}\\
0 & \text{otherwise.}
\end{cases}
$$
In particular, if $\Delta = \cL$ (resp. $\Delta = \{ \, \{ 0 \} \, \}$), 
then $I_\Delta = 0$ (resp. $I_\Delta = \m$) 
and $K[\Delta] = R$ (resp. $K[\Delta] = K$).  
We have  $\dim K[\Delta] = \dim |\Delta| +1$. 
When $R$ is a polynomial ring, $K[\Delta]$ is nothing other than 
the Stanley-Reisner ring of a simplicial complex $\Delta$. 
(If $R$ is simplicial, then $B$ is a simplex and   
$\Delta$ is a simplicial complex.) 
Clearly, $I_\Delta$ and $K[\Delta]$ are squarefree modules. 

\medskip

Let $\Delta, \Sigma \subset \cL$ be order ideals with $\Delta \supset \Sigma$. 
When we consider such a pair, we assume that 
$\Sigma \ne \{ \, \{0 \} \, \}$, but allow the case $\Sigma = \emptyset$. 
Thus, $|\Sigma| = \emptyset$ if and only if 
$\Sigma = \emptyset$. We set $I_\emptyset = R$. 
We always have $I_\Sigma \supset I_\Delta$, and 
$I_{\Delta/\Sigma} := I_\Sigma/I_\Delta$ is a squarefree $R$-module. 
Note that $I_{\Delta/\emptyset} = K[\Delta]$. 
We say a pair $(\Delta, \Sigma)$ is {\it Cohen-Macaulay} 
if so is the module $I_{\Delta/\Sigma}$. 
(The Cohen-Macaulay property of $(\Delta, \Sigma)$ depends on  
$\chara(K)$. But, in this paper, we fix the base field $K$. 
So we omit the phrase ``over $K$".) 
Clearly, the Cohen-Macaulay property of $(\Delta, \Sigma)$ 
depends only on $\Delta \setminus \Sigma$. 
For a topological space $X$, $\Sh(X)$ denotes 
the category of sheaves of $K$-vector spaces on $X$. 
Recall that $Z := |\Delta| \setminus |\Sigma|$ admits {\it Verdier's 
dualizing complex} (over $K$) $\Dc_Z \in D^b(\Sh(Z))$. 
See \cite[V. \S2]{I}. Recall a few results from \cite{Y6}, 
some parts of which have been obtained by Stanley \cite[\S\S4,5]{St2} 
(this does not mean that the next theorem contains 
corresponding results in \cite[\S\S4,5]{St2}, 
since the rings treated there are more general than ours). 

\begin{thm}[{\cite[Theorem~4.9, Proposition~4.10]{Y6}}]\label{topological}
We have the following. 
 \begin{itemize}
\item[(1)]  $K[\Delta]$ is Cohen-Macaulay if and only if 
$\cH^{-i}(\Dc_{|\Delta|}) =0$ and $\rH^i(|\Delta|;K)=0$ for all 
$i \ne \dim |\Delta|$. 
\item[(2)] Let $(\Delta, \Sigma)$ be a pair of order ideals of $\cL$ 
with $\Sigma \ne \emptyset$, and $h$ the embedding map 
from $Z := |\Delta| \setminus |\Sigma|$ to its closure $\overline{Z} 
= |\Delta|$. Then $(\Delta, \Sigma)$ is Cohen-Macaulay if and only if 
${\rm R}^{-i}h_*\Dc_Z =0$  and $H_c^i(Z; K) = 0$ 
for all $i \ne \dim Z$. Here $H_c^i(-)$ denotes 
the cohomology with compact support. 
\end{itemize}
Consequently, either $\Sigma = \emptyset$ or $\Sigma \ne \emptyset$,  
the Cohen-Macaulay property of $(\Delta, \Sigma)$ depends only on 
the pair of topological spaces $(|\Delta|, |\Sigma|)$ 
(or even $(\overline{Z}, Z)$). 
\end{thm}

\begin{rem}
We can check the Cohen-Macaulay property of a pair $(\Delta, \Sigma)$ 
explicitly. Recall the combinatorial 
description of $\Du(I_{\Delta/\Sigma}) = 
\RHom_R(I_{\Delta/\Sigma}, \can) \in D^b(\Sq)$
given in  Lemma~\ref{explicit form}. For each $F \in \cL$, the complex 
$[\Du(I_{\Delta/\Sigma})]_{\bc(F)}$ of $K$-vector spaces is isomorphic to 
the complex $D_F^\bullet$ defined by 
$$D_F^i := \bigoplus_{\substack{G \in \Delta \setminus \Sigma, \ 
G \supset F \\ \dim G = d-i}} K e_G,$$
$$\partial: e_G \longmapsto \sum_{\substack{G' \in \Delta \setminus \Sigma, 
\ G \supset G' \supset F \\ \dim G' = \dim G -1}} \pm \, e_{G'}, $$
where $e_G$ is a basis element, and the sign $\pm$ is given 
by an incidence function of $\cL$. Hence 
$(\Delta, \Sigma)$ is Cohen-Macaulay if and only if 
$H^i(D_F^\bullet)=0$ for all $i \ne d-\dim I_{\Delta/\Sigma}$ 
and all $F \in \Delta$.  
\end{rem}

For $F \in \Delta$, set 
$\delta(F):= \max \{ \, \dim |G| \mid G \supset F, \ G \in \Delta \, \}$. 
For $i \in \NN$ with $i \leq \dim |\Delta|$, we call 
$\Delta^{(i)} := \{ \, F \in \Delta \mid \dim |F| \leq i \, \}$ 
the $i$-{\it skeleton} of $\Delta$, and 
$\Delta^{[i]}:= \{ \, F \in \Delta^{(i)} \mid \delta(F) \geq i \, \}$
the {\it pure $i$-skeleton} of $\Delta$. 
Clearly, these are order ideals of $\cL$ again. 
It is easy to see that 
\begin{equation}\label{pure skeleton}
\Delta^{[i]} \setminus (\Delta^{[i+1]})^{(i)} 
= \{ \, F \in \Delta \mid  \delta(F) = i  \, \}
\end{equation}
for all $i \leq \dim |\Delta| =:r$. 
Here we set $(\Delta^{[r+1]})^{(r)} = \emptyset$. 

For a finitely generated graded $R$-module $M$ and an integer $i$, we set 
$$M_{(i)}:= \{ \, y \in M \mid \dim (Ry) \leq i \, \}.$$
Then $M_{(i)}$ is a submodule of $M$ with $M=M_{(d)} \supset M_{(d-1)} 
\supset \cdots \supset M_{(0)} \supset M_{(-1)} =0$. 
It is easy to see that $M$ is sequentially Cohen-Macaulay 
if and only if $M_{(i)}/M_{(i-1)}$ is Cohen-Macaulay 
(of course, we allow the case $M_{(i)}/M_{(i-1)}=0$) for all $i$. 

Since $K[\Delta]_{(i+1)}/K[\Delta]_{(i)} 
\cong I_{\Delta^{[i]} / (\Delta^{[i+1]})^{(i)} }$, we have the following.

\begin{lem}[c.f. {\cite[Proposition~2.10]{St}, \cite[Remark in p.4]{D}}]
\label{filtration}With the above notation, 
$K[\Delta]$ is sequentially Cohen-Macaulay if and only if the pair
$( \, \Delta^{[i]}, (\Delta^{[i+1]})^{(i)} \,  )$ is 
Cohen-Macaulay for all $i \leq \dim |\Delta|$.  
\end{lem}

When $R$ is a polynomial ring (i.e., $K[\Delta]$ 
is a Stanley-Reisner ring), the next result has been 
obtained by Duval~\cite{D}.  

\begin{thm}[{c.f. Duval \cite[Theorem~3.3]{D}}]\label{duval}
For an order ideal $\Delta \subset \cL$, 
$K[\Delta]$ is sequentially Cohen-Macaulay if and only if  
$K[\Delta^{[i]}]$ is Cohen-Macaulay for all $i \leq \dim |\Delta|$. 
\end{thm}

\begin{proof}
The proof of \cite[Theorem~3.3]{D} for the Stanley-Reisner ring case 
also works here. This proof uses Lemma~\ref{filtration} and 
Lemma~\ref{skelton} below. The latter, which is a well-known result 
in the Stanley-Reisner ring case, also holds in our context. 
\end{proof}

\begin{lem}\label{skelton}
If $K[\Delta]$ is Cohen-Macaulay, then so is $K[\Delta^{(i)}]$ for all 
$i \leq \dim |\Delta|$. 
\end{lem}

\begin{proof}
We may assume that $i+1 = \dim |\Delta| =:r$. 
For each $F \in \Delta$, the canonical module $\omega_{K[F]}$ of 
$K[F] :=R/P_F$ is a squarefree $R$-module with 
$$
[\omega_{K[F]}]_\bc = \begin{cases}
K & \text{if $s(\bc) = F$,} \\
0 & \text{otherwise,}
\end{cases}
$$
for all $\bc \in C$. We have a short exact sequence 
\begin{equation}\label{skeleton}
0 \to \bigoplus_{\substack{F \in \Delta \\ \dim |F| =r}} 
\omega_{K[F]} \to K[\Delta] \to K[\Delta^{(r-1)}] \to 0.
\end{equation}
Since $\omega_{K[F]}$ for $F$ with $\dim |F| = r$ is a Cohen-Macaulay module 
of dimension $r+1$ $(= \dim K[\Delta])$ and $\dim K[\Delta^{(r-1)}] = r$, 
$K[\Delta^{(r-1)}]$ is Cohen-Macaulay.  
\end{proof}

\begin{rem}
In the Stanley-Reisner ring case, Hibi (\cite[\S2]{Hi}) showed that if 
$i < \dim |\Delta|$ then the canonical module of 
$K[\Delta^{(i)}]$ is generated by its degree 0 part. 
This also holds in our context. To see this, we may assume that 
$i+1 = \dim |\Delta| =:r$. Then the long exact sequence of 
$\Ext^i_R(-,\can)$ derived from the sequence \eqref{skeleton} gives 
$$ \bigoplus_{\substack{F \in \Delta \\ \dim |F| =r}} K[F] \to 
\Ext^{d-r}_R(K[\Delta^{(i)}], \can) \to 0 \quad \text{(exact)}.$$
The module of right side, which is isomorphic to the canonical module of 
$K[\Delta^{(i)}]$, is clearly generated by its degree 0 part. 
\end{rem}

Since $|\Delta^{[i]}|$ depends on the cell decomposition 
$|\Delta| = \bigcup_{F \in \Delta} |F|$, Theorem~\ref{duval} does {\it not}  
implies that the sequentially Cohen-Macaulay property of $K[\Delta]$ 
is a topological property of $|\Delta|$. But we can prove this 
using Lemma~\ref{filtration} directly. 

\begin{thm}\label{seqcm}
The sequentially Cohen-Macaulay property of $K[\Delta]$ only depends on 
the topological space $|\Delta|$. 
\end{thm}

\begin{proof}
Set $X:= |\Delta|$. Note that the subspace 
$$\bigcup_{F \in \Delta^{[i]}} |F| \ \, \setminus \ 
\bigcup_{F \in (\Delta^{[i+1]})^{(i)}} |F|$$
of $X$ does not depend on the particular cell decomposition of $X$ by 
\eqref{pure skeleton}. In fact, it coincides with 
$\{ \, x \in X \mid \dim_x X = i \, \}$, where the dimension 
$\dim_x X$ of $X$ at $x$ is the one defend in \cite[III, Definition~9.10]{I}. 
So the assertion follows from Theorem~\ref{topological} 
and Lemma~\ref{filtration}. 
\end{proof}

Even in the Stanley-Reisner ring case, it is well-known that 
the Gorenstein property of $K[\Delta]$ is not a topological property of 
$|\Delta|$ (i.e., depends on the simplicial decomposition). 
In the normal semigroup ring case, we have another problem. 
A normal semigroup ring $R$ is always Cohen-Macaulay, 
but not necessarily Gorenstein, even if it is simplicial. 
Note that $K[\Delta] = R$ if $\Delta = \cL$. 
So the Gorenstein property can not be determined by the poset 
structure of $\Delta$ (i.e., sensitive to the semigroup $C$). 
But we can prove that the {\it Gorenstein*} property is topological. 
Recently, Ichim and R\"omer also studied the Gorenstein (or Gorenstein*)
property of a {\it toric face ring}, which is a notion containing our  
$K[\Delta]$. 
Their \cite[Corollary~6.9]{IR} is closely related to 
Theorem~\ref{Gor*} below.

\begin{dfn}[{Stanley, \cite[p.67]{St}}]
We say $K[\Delta]$ is {\it Gorenstein*}, if it is Gorenstein 
and the canonical module $\omega_{K[\Delta]}$ of $K[\Delta]$ 
is generated by its degree 0 part (equivalently,  
$\omega_{K[\Delta]} \cong K[\Delta]$ as 
graded modules).  
\end{dfn}

\begin{prop}\label{Gor*}
Set $X := |\Delta|$, and let $\Ddel \in D^b(\Sh(X))$ 
be the dualizing complex of $X$. 
Assume that $r:= \dim X \geq 1$ and $K[\Delta]$ is Cohen-Macaulay. 
Then $K[\Delta]$ is Gorenstein* if and only if 
$\cH^{-r}(\Ddel)_x =  K$ for all $x \in X$ and $H^r(X; K) \ne 0$. 
Thus the Gorenstein* property of $K[\Delta]$ 
only depends on the topological space $X$. 
\end{prop}

Under the assumption of Proposition~\ref{Gor*}, we have 
$\cH^{-r}(\Ddel)[r] \cong \Ddel$ in $D^b(\Sh(X))$. 
When $X$ is a manifold with or without boundary 
(but  $K[\Delta]$ is not necessarily Cohen-Macaulay), 
$\cH^{-r}(\Ddel)$ is called the {\it orientation sheaf} 
of $X$ (over $K$). Hence, when $X$ is a manifold and $K[\Delta]$ 
is Cohen-Macaulay, $K[\Delta]$ is Gorenstein* if and only if 
$X$ is an orientable manifold over $K$
by \cite[Theorem~4.2]{Y6} (quoted as Theorem~\ref{Verd} below). 
In particular, if $X$ is homeomorphic to an $r$-dimensional sphere, 
then $K[\Delta]$ is Gorenstein*. 

\begin{proof}
The last statement easily follows from the first one. 
In fact, the assumption that $r \geq 1$ is not essential. 
If $r = -1$ (i.e., $X = \emptyset$), then $K[\Delta]$ 
is Gorenstein*. When $r=0$, $K[\Delta]$ is Gorenstein* if and only if 
$X$ consists of two points. On the other hand, we have seen that 
the Cohen-Macaulay property of $K[\Delta]$ is a topological property of $X$. 
So it suffices to prove the first statement. 
Recall that $\omega_{K[\Delta]}$ is a squarefree module with 
$\Ass(K[\Delta]) = \Ass(\omega_{K[\Delta]})$. By \cite[Theorem~4.2]{Y6}, 
$\cH^{-r}(\Ddel)_x = K$ for all $x \in X$ if and only if 
$[\omega_{K[\Delta]}]_{\bc(F)} = K$ for all $F \in \Delta \setminus 
\{ 0 \}$. Hence, under the assumption that 
$\cH^{-r}(\Ddel)_x = K$ for all $x \in X$, we have 
$\omega_{K[\Delta]} \cong K[\Delta]$ as graded modules, if and only if 
$\omega_{K[\Delta]}$ is generated by its degree 0 part, if and only if 
$[\omega_{K[\Delta]}]_0 \ne 0$ (equivalently, 
$[\omega_{K[\Delta]}]_0 = K$).  
By \cite[Theorem~3.3]{Y6} (quoted as Theorem~\ref{Hoch} below) 
and the local duality, we have 
$H^r(X;K) \cong [H_\m^{r+1}(K[\Delta])]_0 \cong 
[\omega_{K[\Delta]}]_0^\vee$. So we are done. 
\end{proof}

\begin{exmp}
In Proposition~\ref{Gor*}, that    
$\cH^{-r}(\Ddel)_x =K$ for all $x \in X$ is not a sufficient condition 
for $K[\Delta]$ to be Gorenstein*. In fact, if $X$ is a manifold 
without boundary, then $\cH^{-r}(\Ddel)_x = K$ for all $x \in X$. 
But, of course, $X$ need not be orientable. 
For example, if $X$ is homeomorphic to a real projective plane 
and $\chara(K) \ne 2$, then $K[\Delta]$ is Cohen-Macaulay but not 
Gorenstein*. 
\end{exmp}

\section{Constructible sheaves associated with squarefree modules}
Throughout this section, $R$ is normal, and  
$\Delta$ is an order ideal of $\cL$ with $X := |\Delta|$. 

If $F \in \Delta$, then $U_F := 
X \cap (\bigcup_{F' \supset F} |F'|)$ is an open set of $X$. 
Note that $\{ \, U_F \mid \{ 0 \} \ne F \in \Delta \, \}$ 
is an open covering of $X$. 
If $M$ is a squarefree $K[\Delta]$-module (i.e., $M \in \Sq$ and 
$\ann(M) \supset I_\Delta$), then we can construct a sheaf 
$M^+$ on $X$ as in \cite{Y6} (see also Remark~\ref{sheaf} (1) below). 
More precisely, the assignment 
$\Gamma(U_F, M^+) = M_{\bc(F)}$ for each $\{ 0  \} \ne F \in \cL$ and 
the map $\varphi_{F,F'}^M: M^+(U_{F'}) = M_{\bc(F')} \to M_{\bc(F)} = 
\Gamma(U_F,M^+)$ for $\{  0 \} \ne F,F' \in \Delta$ with 
$F \supset F'$ (equivalently, $U_{F'} \supset U_{F}$) actually defines 
a sheaf. It is easy to see that $M^+$ is a {\it constructible sheaf} 
with respect to the cell decomposition $X=\bigcup_{F \in \Delta} |F|$. 
In fact, for all $\{ 0 \} \ne F \in \Delta$, 
the restriction $M^+|_{|F|}$ of $M^+$ to $|F| \subset X$ 
is a constant sheaf with coefficients in $M_{\bc(F)}$.  
Thus, if  $p \in |F| \subset X$, then the stalk $(M^+)_p$ at $p$ is 
isomorphic to $M_{\bc(F)}$. If $M$ is a squarefree module, 
then so is the submodule $M_{>0} := \bigoplus_{0 \ne \bc \in C}M_\bc$. 
For squarefree $K[\Delta]$-modules $M$ and $N$, it is easy to see that 
$M^+ \cong N^+$ if and only if $M_{> 0} \cong N_{> 0}$. 
In other words, $M_0$ is ``irrelevant" to $M^+$. 

Let $\Sq(\Delta)$ be the category of squarefree $K[\Delta]$-modules 
(hence $\Sq(\Delta)$ is a full subcategory of $\Sq$).  
The above operation gives the exact functor $(-)^+: \Sq(\Delta) \to \Sh(X)$. 

\begin{exmp}
Note that $K[\Delta]$ itself is a squarefree $K[\Delta]$-module, and we have
$K[\Delta]^+ \cong \const$, where $\const$ is the 
constant sheaf on $X:= |\Delta|$ with coefficients in $K$.  
For a face $\{ 0 \} \ne F \in \Delta$, we set 
$K[F] := K[\Delta]/P_F$. 
Then we have $K[F]^+ \cong j_* \const_{|F|} \in \Sh(X)$, where 
$j: |F| \hookrightarrow X$ is the embedding map and 
$\const_{|F|}$ is the constant sheaf on $|F|$ 
(note that $j_* \const_{|F|}$ is essentially the constant sheaf on the 
closure $\overline{|F|}$ of $|F|$, not on $|F|$ itself). 
Recall that the canonical module $\can$ is a squarefree $R$-module, and 
$R=K[\Delta]$ if $\Delta = \cL$. 
Then $\can^+ \cong h_! \const_{B^\circ} \in \Sh(B)$, 
where $B^\circ$ is the interior of $B$, $h: B^\circ \hookrightarrow B$ 
is the embedding map, and $\const_{B^\circ}$ is the constant sheaf on 
$B^\circ$. 
It is noteworthy that $h_! \const_{B^\circ}$ is the {\it orientation sheaf} 
of $B$ (over $K$). 
This is a key point for Theorem~\ref{Verd} below. 
\end{exmp}

\begin{thm}[{\cite[Theorem~3.3]{Y6}}]\label{Hoch}
For $M \in \Sq(\Delta)$, we have an isomorphism 
$$H^i(X; M^+) \cong [H_\m^{i+1}(M)]_0 \quad  \text{for all $i \geq 1$},$$ 
and an exact sequence 
\begin{equation}\label{lower dimension}
0 \to [H_\m^0(M)]_0 \to M_0 \to H^0( X; M^+) \to [H_\m^1(M)]_0 \to 0.
\end{equation}
In particular, we have 
$[H_\m^{i+1}(k[\Delta])]_0 \cong \rH^i(X ; K)$ for all $i \geq 0$, 
where $\rH^i(X ; K)$ denotes the $i^{ th}$ reduced cohomology 
of $X$ with coefficients in $K$. 
\end{thm}

\begin{rem}\label{sheaf}
(1) In \cite{Y6}, we construct the sheaf $M^+$ regarding $M$ 
as an $R$-module (not a $K[\Delta]$-module). 
Thus, $M^+$ is always a sheaf on $B$ there. But this is not a problem. 
In fact, if $\cF \in \Sh(B)$ is the sheaf constructed from $M \in \Sq(\Delta)$ 
in the style of \cite{Y6}, and $M^+ \in \Sh(X)$ is the sheaf constructed 
in the style of the present paper, 
then we have $\cF \cong i_* M^+$. 
Here $i : X \hookrightarrow B$ is the embedding map. 
Since $H^i(X; M^+) \cong H^i(B; \cF)$ for all $i$, 
\cite[Theorem~3.3]{Y6} implies Theorem~\ref{Hoch}.

(2) When $R$ is a polynomial ring (i.e., $K[\Delta]$ is 
a Stanley-Reisner ring), the last statement of Theorem~\ref{Hoch}  
is (a part of) a famous formula of Hochster 
(\cite[Theorem~II.4.1]{St}). 
We can also compute $[H_\m^i(M)]_\ba$ for any $\ba \in \Z^d$ 
in terms of the sheaf $M^+$ (see \cite[Theorem~3.5]{Y6}), 
while we do not use this result here. 

(3) For $y \in M_0$, the sections 
$\varphi_{F, \{ 0 \} }^M (y) \in M_{\bc(F)} = \Gamma(U_F; M^+)$ 
for all $\{ 0\} \ne F \in \cL$ define the unique global section 
$\tilde{y} \in \Gamma(X ;M^+)$. The map $M_0 \to H^0( X, M^+)$ in 
the sequence \eqref{lower dimension} 
is given by $y \mapsto \tilde{y}$. 
\end{rem}

\begin{dfn}
A squarefree module $M$ is {\it regular}, if $[H_\m^0(M)]_0 =
[H_\m^1(M)]_0 = 0$. 
\end{dfn}

\begin{lem}\label{regular}
Let $M$ be a squarefree $K[\Delta]$-module. There is a unique 
(up to isomorphism) squarefree $K[\Delta]$-module $\overline{M}$ which 
is regular and $\overline{M}^+ \cong M^+$.  
\end{lem}

\begin{proof}
We put $\overline{M}_{>0} = M_{>0}$ and $\overline{M}_0 = \Gamma(X; M^+)$. 
It suffices to define $\varphi_{F,\{ 0\}}^{\overline{M}}: \overline{M}_0 \to 
\overline{M}_{\bc(F)}$ for each $\{ 0\} \ne F \in \cL$. 
Under the identification $\overline{M}_0 = \Gamma(X; M^+)$ and 
$\overline{M}_{\bc(F)} = M_{\bc(F)} = \Gamma(U_F; M^+)$, 
$\varphi_{F,\{ 0\}}^{\overline{M}}$ 
is given by the restriction map $\Gamma(X ; M^+) \to \Gamma(U_F ; M^+)$. 
Then $\overline{M}$ satisfies the expected condition by the sequence 
\eqref{lower dimension}. 
\end{proof}

\begin{rem}
If $X$ is connected, then we have $K[\Delta] \cong \overline{K[\Delta]}$ 
in the notation of Lemma~\ref{regular}. 
But, if $X$ is not connected, then 
$\overline{K[\Delta]}_0 = K^n$, where $n$ is the number of the 
connected components  of $X$. 
\end{rem}

Set $r:=\dim X$. Consider the complex 
$$\wdel: 0 \to \Wde^{-r} \to \Wde^{-r+1} \to 
\cdots \to \Wde^0 \to \Wde^1 \to 0,$$ 
$$\Wde^i := \bigoplus_{\substack{F \in \Delta \\ \dim |F| = -i}} K[F]$$
of squarefree $K[\Delta]$-modules. The translated complex 
$\wdel[1]$ of $\wdel$ is (quasi-isomorphic to) the {\it normalized} 
dualizing complex of $K[\Delta]$. Hence we have 
$$\Ext_{K[\Delta]}^i(\M, \wdel) \cong \Ext_{K[\Delta]}^{i-1}(\M, \wdel[1])
\cong H_\m^{-i+1}(\M)^\vee$$ for all $\M \in D^b(\Sq(\Delta))$ 
and all $i \in \Z$.  
But here we prefer $\wdel$ to $\wdel[1]$, since 
$(\wdel)^+ \in D^b(\Sh(X))$ is quasi-isomorphic to 
the dualizing complex $\Ddel$ of $X$ as shown in \cite[corollary~4.3]{Y6}.  

If $\M \in D^b(\Sq(\Delta))$, then we have $\RHom_{K[\Delta]}(\M, \wdel) 
\in D^b_{\Sq(\Delta)}(\GrR) \cong D^b(\Sq(\Delta))$. 
So we can define $\RHom_{K[\Delta]}(\M, \wdel)^+ \in D^b(\Sh(X))$. 
Moreover, the following result holds. 

\begin{thm}[{\cite[Theorem~4.2]{Y6}}]\label{Verd}
For $\M \in D^b(\Sq(\Delta))$, we have 
$$\RHom_{K[\Delta]}(\M, \wdel)^+ 
\cong {\rm R}\cHom ((\M)^+, \Ddel)$$ in $D^b(\Sh(X))$. 
In particular, $\Ext^i_{K[\Delta]}(\M, \wdel)^+ \cong 
\cExt^i((\M)^+, \Ddel)$. 
\end{thm}

The next two results easily follow from Theorem~\ref{Verd} and 
the local duality theorem for $K[\Delta]$. 

\begin{prop}\label{non-vanishing}
If $M$ is a squarefree $K[\Delta]$-module with $\dim M \geq 1$, 
then $$ - \min \{ \, i \mid \cExt^i(M^+, \Ddel) \ne 0 \, \} 
= \dim \Supp(M^+) = \dim M -1.$$
Here $\Supp(M^+) := \{ \, x \in X \mid (M^+)_x \ne 0 \, \}$.  
\end{prop}

\begin{thm}\label{CM sheaf}
Let $M$ be a squarefree $K[\Delta]$-module which is regular. 
Set $t: = \dim \Supp(M^+) = \dim M-1$. 
Then $M$ is Buchsbaum if and only if 
$\cExt^i(M^+, \Ddel) = 0$ for all $i \ne -t$. 
Similarly, $M$ is Cohen-Macaulay if and only if $\cExt^i(M^+, \Ddel) = 0$ 
for all $i \ne -t$ and $H^i(X; M^+)=0$ for all $i \ne 0, t$. 
In particular, the Buchsbaum property and the Cohen-Macaulay property of 
$M$ depend only on the sheaf $M^+$. 
\end{thm}

\begin{rem}
In Theorem~\ref{CM sheaf} (and Theorem~\ref{seqcm module} below), 
even $X$ is somewhat superfluous, and 
the {\it closure} of $\Supp(M^+)$ is essential. 
If we set $T:=\{ \, F \in \Delta \mid M_{\bc(F)} \ne 0 \, \}$, 
then $\bigcup_{F \in T} |F| = \Supp(M^+)$. 
Let $\Sigma$ be the order ideal of $\cL$ generated by $T$. 
Then $Y:=\bigcup_{F \in \Sigma} |F|$ coincides with 
the closure of $\Supp(M^+)$. 
Note that $\Sigma \subset \Delta$ and $M$ can be regarded as a squarefree 
$K[\Sigma]$-module. Let $\cF \in \Sh(Y)$ be the sheaf associated with 
$M$ as a $K[\Sigma]$-module. Then $\cF \cong M^+|_Y$, 
$H^i(Y; \cF) \cong H^i(X; M^+)$, and 
$\cExt^i(\cF, \Dsig) \cong \cExt^i(M^+, \Ddel)$ for all $i$. 
Of course, the Cohen-Macaulay (resp. Buchsbaum) property of $M$ 
does  not depend on whether we regard $M$ as a $K[\Delta]$-module 
or a $K[\Sigma]$-module.  
\end{rem}

Let $\cF \in \Sh(X)$, and set $\cG := {\rm R}\cHom(\cF, \Ddel) \in 
D^b(\Sh(X))$. Recall that $\cExt^i(\cF, \Ddel) = \cH^i(\cG)$. 
More precisely, $\cExt^i(\cF, \Ddel)$ 
is the sheaf associated with the presheaf defined by 
$U \mapsto H^i(\Gamma(U; \cG))$ for an open subset $U$ of $X$. 
Hence the element of $\Ext^i(\cF, \Ddel) = H^i(\Gamma(X; \cG))$ 
gives a global section of $\cExt^i(\cF, \Ddel)$, that is, 
we have a natural map 
$\Ext^i(\cF, \Ddel) \to \Gamma(X; \cExt^i(\cF, \Ddel))$. 

\begin{lem}\label{wdel}
If $M \in \Sq(\Delta)$, then $[\Ext^i_{K[\Delta]}(M, \wdel)]_0 \cong 
\Ext^i(M^+, \Ddel)$ for all $i < 0$. Via this isomorphism, 
the natural map $\Ext^i(M^+, \Ddel) \to \Gamma(X; \cExt^i(M^+, \Ddel))$ 
coincides with the middle map $[\Ext^i_{K[\Delta]}(M, \wdel)]_0 \to 
\Gamma(X; \Ext^i_{K[\Delta]}(M, \wdel)^+)$ of the complex \eqref{lower dimension} 
in Theorem~\ref{Hoch}. 
\end{lem}

\begin{proof} 
Set $\N := \RHom_{K[\Delta]}(M,\wdel) \in D^b(\Sq(\Delta))$. 
By the same argument as Lemma~\ref{explicit form}, 
we may assume that $$N^i = \bigoplus_{\substack{F \in \Delta 
\\ \dim |F| = -i}}(M_{\bc(F)})^* \otimes_K K[F].$$ 
Thus $(N^i)^+$ is a direct sum of copies of the sheaf $K[F]^+$ for 
various $F \in \Delta$ with $\dim |F| = -i$. 
While $K[F]^+$ is not an injective object in $\Sh(X)$, 
it is a constant sheaf over the closed ball $|F|$ and 
$H^i(X;K[F]^+)=H^i(|F|;K)=0$ for all $i>0$. From this fact and that 
$(\N)^+ \cong \RcHom (M^+, \Ddel)$,  
$\Ext^i(M^+, \Ddel)$ is the $i^{\rm th}$ cohomology of the complex 
$\Gamma(X; (\N)^+)$. On the other hand, 
$\Gamma(X; (N^i)^+) = [N^i]_0$ for all $i \leq 0$. 
Hence $\Ext^i(M^+, \Ddel) = [\Ext_{K[\Delta]}^i(M,\wdel)]_0$ for all 
$i < 0$. Since $H^i(\N)^+ \cong \cExt^i (M^+, \wdel)$, 
the assertion can be checked easily. 
\end{proof}

\begin{thm}\label{seqcm module}
Assume that $M \in \Sq(\Delta)$ is regular. Then $M$ is sequentially 
Cohen-Macaulay if and only if the following conditions are satisfied.
\begin{itemize}
\item[(a)] $\cExt^i(\cExt^j(M^+, \Ddel), \Ddel) =0$ for all 
$i,j \in \Z$ with $i \ne j$ and $j < 0$. 
\item[(b)] $H^i( X ;\cExt^j(M^+, \Ddel)) =0$ 
for all $i \ne 0, -j$ and all $j<0$.
\item[(c)] The natural map $\Ext^i(M^+, \Ddel) \to \Gamma(X; 
\cExt^i(M^+, \Ddel))$ is bijective for all $i < 0$.  
\end{itemize}
In particular, the sequentially Cohen-Macaulay property of $M$ 
only depends on the sheaf $M^+$. 
\end{thm}

\begin{proof}
By argument similar to \cite[Theorem~III.2.11]{St}, we can see that 
$M$ is sequentially Cohen-Macaulay, if and only if 
$\Ext_{K[\Delta]}^i(M, \wdel)$ is either the 0 module  
or a Cohen-Macaulay module of dimension $1-i$ for all $i$. 
Since  $\Ext_{K[\Delta]}^i(M, \wdel) = 0$ for $i > 1$ and 
$\Ext_{K[\Delta]}^1(M,\wdel)$ is an artinian module, we do not have to check 
$\Ext_{K[\Delta]}^i(M,\wdel)$ for $i \geq 1$. Moreover, the case when 
$i = 0$ is also unnecessary in our situation. 
In fact, since $M$ is regular, we have 
$[H_\m^1(M)^\vee]_0 = [\Ext_{K[\Delta]}^0(M,\wdel)]_0=0$.  
From this fact and that $\dim \Ext_{K[\Delta]}^0(M,\wdel) \leq 1$, 
$\Ext_{K[\Delta]}^0(M,\wdel)$ is Cohen-Macaulay. 

If the condition (c) is satisfied, $\Ext_{K[\Delta]}^i(M,\wdel)$ 
is regular by Lemma~\ref{wdel}. 
Conversely, if $M$ is sequentially Cohen-Macaulay, then
$\Ext_{K[\Delta]}^i(M,\wdel)$ must be regular for $i < 0$ 
and (c) is satisfied.  

Since $\Ext_{K[\Delta]}^i(M, \wdel)^+ \cong \cExt^i(M^+, \Ddel)$, 
the assertion follows from the above observation and 
Theorem \ref{CM sheaf} (and Proposition \ref{non-vanishing}).
\end{proof}

\begin{rem}
If $X$ is not connected, then $K[\Delta]$ is not regular 
as a squarefree module. 
So Theorem~\ref{seqcm module} does not imply Theorem~\ref{seqcm} 
directly. But, by the following observation, Theorem~\ref{seqcm} 
can be reduced to Theorem~\ref{seqcm module}. 
\begin{itemize}
\item[(1)] 0 dimensional components of $X$ 
(i.e., 1 dimensional components of $K[\Delta]$) are irrelevant 
to the sequentially Cohen-Macaulay property of $K[\Delta]$. 
So we can remove them.  
\item[(2)] If $X$ does not have a 0 dimensional component and 
$K[\Delta]$ is sequentially Cohen-Macaulay, then $X$ is connected. 
\end{itemize} 
\end{rem}

\section{Ideals whose radicals are monomial ideals}
In this brief section, we generalize a result of 
Herzog, Takayama and Terai~\cite{HTT}. 
Let $I \subset R$ be a (non-monomial) ideal. 
Even if $R/I$ is Cohen-Macaulay,  
$R/\sqrt{I}$ is {\it not} Cohen-Macaulay in general. See the introduction of 
\cite{HTT} for an explicit example. But the next theorem states that 
if $\sqrt{I}$ is a monomial ideal then such an example does not exist. 
When $R$ is a polynomial ring, this result was obtained in \cite{HTT}. 

\begin{thm}[c.f. {\cite[Theorem~2.6]{HTT}}] 
Assume that $R$ is normal. 
Let $I$ be a (not necessarily graded) ideal whose 
radical $\sqrt{I}$ is a monomial ideal. 
Then we have the following 

(1)If $R/I$ is Cohen-Macaulay (more generally, the localization $(R/I)_\m$ 
is Cohen-Macaulay), then $R/\sqrt{I}$ is also. 

(2) If  the localization $(R/I)_\m$ is generalized Cohen-Macaulay, 
then $R/\sqrt{I}$ is Buchsbaum, in particular, it is 
generalized Cohen-Macaulay again. 
\end{thm}

The idea of the proof is same as the one given in \cite[Remark~2.7]{HTT}.

\begin{proof}
Let me introduce the notation and facts used throughout this proof. 
Set $p:= \dim (R/I) = \dim (R/\sqrt{I})$.
Since $\sqrt{I}$ is a monomial ideal, $\Ext^i_R(R/\sqrt{I}, \can)$ 
has a natural $\Z^d$-grading. In particular, 
$\Ext^i_R(R/\sqrt{I}, \can) \otimes_R R_\m = 0$ implies 
$\Ext^i_R(R/\sqrt{I}, \can) = 0$. 
Similarly, if $\Ext^i_R(R/\sqrt{I}, \can) \otimes_R R_\m$ has finite length, 
then so does $\Ext^i_R(R/\sqrt{I}, \can)$. 

(1) It suffices to show that $\Ext_R^i(R/\sqrt{I}, \can) = 0$ 
for all $i \ne d- p$.  Recall that the natural map 
$$\Ext_R^i(R/\sqrt{I}, \can) \to H_{\sqrt{I}}^i(\can)$$ 
factors through the map 
$f: \Ext_R^i(R/\sqrt{I}, \can) \to \Ext_R^i(R/I, \can)$ 
induced by the surjection $R/I \to R/\sqrt{I}$. 
But the map $\Ext_R^i(R/\sqrt{I}, \can) \to H_{\sqrt{I}}^i(\can)$ 
is an injection by \cite[Theorem~5.9]{Y2}. 
So the map $f$ and its localization $f \otimes_R R_\m$  
are injective. Since $(R/I)_\m$ is Cohen-Macaulay, we have 
$\Ext_R^i(R/I, \can) \otimes_R R_\m = 0$ for all $i \ne d- p$. 
So we are done. 

(2) By the assumption, 
$\Ext_R^i(R/I, \can) \otimes_R R_\m \cong H_\m^{d-i}((R/I)_\m)^\vee$ 
has finite length for all $i \ne d- p$. 
Recall that the map $f: \Ext_R^i(R/\sqrt{I}, \can) \to \Ext_R^i(R/I, \can)$ 
and the localization $f \otimes R_\m$ are injective. 
Hence if $\Ext_R^i(R/\sqrt{I}, \can)$ 
does not have finite length, then $\Ext_R^i(R/I, \can) \otimes_R R_\m$ also. 
So $\Ext_R^i(R/\sqrt{I}, \can)$ has finite length for all $i \ne d- p$. 
By \cite[Corollary~4.6]{Y6}, $R/\sqrt{I}$ is Buchsbaum. 
\end{proof}

\section{Local cohomology modules of finite length}
Let $S=K[x_1, \ldots, x_d]$ be a polynomial ring, and $I \subset S$ a monomial 
ideal. Recently, Takayama~\cite{T} gave an interesting observation that 
the range  
$$\{ \, \ba \in \Z^d \mid [H_\m^i(S/I)]_\ba \ne 0 \, \}$$ 
is controlled by the degrees of minimal generators of $I$ 
(especially, when $H_\m^i(S/I)$ has finite length).  
For this result, he used a combinatorial formula on  $[H_\m^i(S/I)]_\ba$. 
But there is an easy and conceptual proof, and we can generalize his result 
to finitely generated $\NN^d$-graded $S$-modules. We regard 
$\NN^d$ as a partially ordered set with $\ba \leq \bb$ 
$\stackrel{\text{def}}{\Longleftrightarrow}$ $a_i \leq b_i$  for all $i$. 
Here $\ba = (a_1, \ldots, a_d)$ and $\bb = (b_1, \ldots, b_d)$. 

\begin{dfn}[{Miller \cite{Mil0}}]
Let $\ba \in \NN^d$. We say an $S$-module $M$ is  
{\it positively $\ba$ determined} (``$\ba$-p.d.", for short)
if $M$ is $\NN^d$-graded, finitely generated, and the multiplication map 
$M_\bb \ni y \mapsto x_i y \in M_{\bb + \be_i}$ is bijective 
for all $\bb \in \NN^d$ and all $1 \leq i \leq d$ with $b_i \geq a_i$. 
Here $\be_i = (0, \ldots, 1, \ldots, 0) \in \NN^d$ is the vector with 1 at 
the $i^{\rm th}$ position. 
\end{dfn}

\begin{rem}
(1) Any finitely generated $\NN^d$-graded $S$-module is an 
$\ba$-p.d.module for sufficiently large $\ba \in \NN^d$. 

(2) Set $\b1 := (1, \ldots, 1) \in \Z^d$.
Then $\b1$-p.d.modules are nothing other than squarefree modules 
(recall Definition~\ref{sq}) over $S$. 

(3) Let $I = (\bx^{\bb_1}, \ldots, \bx^{\bb_r})$ be a monomial ideal in $S$ 
with $\bb_j= (b^{j}_1, \ldots, b^{j}_d)$. 
Set $a_i := \max \{ \, b^{j}_i \mid 1 \leq j \leq r \, \}$ 
for each $1 \leq i \leq d$, and set $\ba := (a_1, \ldots, a_d)$ 
(in other words, $\bx^\ba= \lcm(\bx^{\bb_1}, \ldots, \bx^{\bb_r})$). 
Then $I$ and $S/I$ are $\ba$-p.d.modules. 
\end{rem}

As shown in \cite{Mil0}, $\ba$-p.d.modules enjoy many interesting properties.  
Some of them are used in the proof of the next result.

\begin{prop}[c.f. {\cite[Corollary~2]{T}}]\label{Takayama}
Let $M$ be an $S$-module which is $\ba$-p.d. 
If $[H_\m^i(M)]_\bb \ne 0$ for some $i$, then $\bb \leq \ba - \b1$. 
\end{prop}

\begin{proof}
Let $\P$ be a $\Z^d$-graded minimal free resolution of $M$. 
Then each $P^i$ is an $\ba$-p.d.module again. That is, 
if $S(-\bb)$ appears in $P^i$ as a direct summand, then 
$\00 \leq \bb \leq \ba$. Let $\wS = S(-\b1)$ be the canonical module of $S$. 
If $S(-\bb)$ is $\ba$-p.d., then 
$$\Hom_S(S(-\bb), \wS)(-\ba+\b1) \cong 
\Hom_S(S(-\bb), S(-\ba)) \cong S(-\ba+\bb)$$ 
is also. 
Hence $[\Hom_S(\P,\wS)](-\ba+\b1)$ is a complex of $\ba$-p.d.modules, and 
its $i^{\rm th}$ cohomology $[\Ext^i_S(M,\wS)](-\ba+\b1)$ 
is also $\ba$-p.d. Since $H_\m^{d-i}(M)_\bb \cong 
[\Ext^i_S(M,\wS)]^\vee_{-\bb}$ by the local duality, we are done. 
\end{proof}

Since $[\Ext^i_S(M,\wS)](-\ba+\b1)$ is an $\ba$-p.d.module 
as shown in the proof of Proposition~\ref{Takayama}, we have the following.

\begin{prop}[c.f. {\cite[Proposition~1]{T}}]\label{FLC}
Let $M$ be a finitely generated $\NN^d$-graded $S$-module. 
Then $H_\m^i(M)$ has finite length if and only if  $H_\m^i(M)$  
is $\NN^d$-graded. 
\end{prop}

To extend Proposition~\ref{FLC} to semigroup rings, 
we have to introduce the notion $\supp(u)$ for $u \in \R^d$. 
There are (the ``defining equations" of) hyperplanes 
$h_1, \ldots, h_t \in (\R^d)^*$ such that 
$\bP \, (= \R_{\geq 0}C) 
= \{ \, u \in \R^d \mid \text{$h_i(u) \geq 0$ for all $i$ }\, \}$. 
We may assume that $h_1, \ldots, h_t$ form a minimal system 
defining $\bP$ (equivalently, the number of $d-1$ dimensional faces 
of $\bP$ is $t$). 
For $u \in \R^d$, set $$\supp(u) := \{ \, i \mid  h_i(u) > 0 \, \} 
\subset \{1, \ldots, t \}.$$ 
For $u, v \in \bP$, $\supp(u)=\supp(v)$ if and only if 
$s(u)=s(v)$. 

So let be $h_1, \ldots, h_t \in (\R^d)^*$ the ``defining equations" 
of the cone $\bP \subset \R^d$. Recall that,  for $u \in \R^d$, we have 
$\supp(u) = \{ \, i \mid  h_i(u) > 0 \, \} \subset \{1, \ldots, t \}$. 
Set $\Cb := \Z^d \cap \bP$. Note that $R$ is normal if and only if $C = \Cb$. 
It is easy to see that, for $\ba \in \Z^d$, 
$\supp(\ba) = \emptyset$ if and only if $\ba \in -\Cb$.  
We say $M \in \GrR$ is $\Cb$-graded, if $M_\ba=0$ for all $\ba \not \in  \Cb$. 
Clearly, a $C$-graded module is always $\Cb$-graded, and the converse is true 
if $R$ is normal. We also set 
$$\Psi_C := \{ \, \ba \in \Z^d \mid \text{$\supp(-\ba) \supset \supp(\bc)$ 
for some $0 \ne \bc \in C$} \, \}.$$

\begin{thm}\label{main}
Assume that $M \in \grR$ is $\Cb$-graded. Then $H_\m^i(M)$ has finite length 
if and only if $[H_\m^i(M)]_\ba = 0$ for all $\ba \in \Psi_C$. 
\end{thm}

\begin{cor}\label{simplicial}
Assume that $R$ is simplicial and $M \in \grR$ is a $\Cb$-graded module.  
Then $H_\m^i(M)$ has finite length if and only if $H_\m^i(M)$ is 
$\Cb$-graded. 
\end{cor}

\begin{proof}
Since $R$ is simplicial, we have 
$$\Psi_C = \{ \ba \in \Z^d \mid \supp(- \ba) \ne \emptyset \}  
= \{ \, \ba \in \Z^d \mid -\ba \not \in  -\Cb \, \} = \Z^d \setminus \Cb.$$ 
So the assertion follows from Theorem~\ref{main}.  
\end{proof}

\begin{exmp}
If $C$ is not simplicial, Corollary~\ref{simplicial} does not hold.  
Let $R$ be the ring given in Example~\ref{square} (1), and set  
$M = R/(x^3, x^2y^2, y^3)$. Then computation by {\it Macaulay 2} 
shows that $H_\m^1(M)$ has finite length 
and $H_\m^1(M)_{(-2,0,0)} \ne 0$. 
\end{exmp}

To prove the theorem, we consider the following condition for 
a module $N \in \GrR$. 
\begin{itemize}
\item[($\star$)] If $\ba \in \Z^d$ and $0 \ne \bc \in C$ satisfy 
$\supp(-\ba ) \supset \supp(\bc)$, then the multiplication map 
$N_{\ba - \bc} \ni y \mapsto \bx^\bc y \in N_{\ba}$ is bijective. 
\end{itemize}

\begin{lem}\label{artinian}
Assume that a $\Z^d$graded $R$-module $N$ has finite length and 
satisfies the condition $(\star)$. Then $N_\ba = 0$ for all $\ba \in \Psi_C$. 
\end{lem}

\begin{proof}
If $\ba \in \Psi_C$, then we can take $0 \ne \bc \in C$ such that  
$\supp(-\ba ) \supset \supp(\bc)$. Since $\supp(-\ba+n\bc) = \supp(\ba)$ 
for all $n \geq 0$, we have 
$N_\ba \cong N_{\ba - \bc} \cong N_{\ba - 2\bc} \cong 
\cdots$ by the condition $(\star)$. So $N_\ba$ must be 0. 
\end{proof}

\begin{lem}\label{5 lemma}
The full subcategory of $\GrR$ consisting of modules 
satisfying $(\star)$ is closed under kernel and cokernels.
\end{lem}

\begin{proof}
Follows from the five-lemma. 
\end{proof}

\noindent{\it Proof of Theorem~\ref{main}.}
(Necessity) \, 
By Lemma~\ref{artinian}, it suffices to prove that $H_\m^i(M)$ satisfies the 
condition $(\star)$. To show this, we use the ``\v{Cech} complex" 
$\L$ of \cite[\S6.2]{BH} (see also \cite[Chapter 13]{MS}. 
There this complex is called 
``Ishida's complex"), 
which satisfies $H_\m^i(M) \cong 
H^i(M \otimes_R \L)$ for all $i$. For a face $F \in \cL$, 
$R_F$ denotes the localization of $R$ at the multiplicatively closed set 
$\{ \, \bx^\bc \mid \bc \in F \cap C \, \}$. Then we have 
$L^i = \bigoplus_{\dim F = i} R_F$  for each $i$. 
 
By Lemma~\ref{5 lemma}, it suffices to prove that $M_F:=M \otimes_R R_F$ 
satisfies $(\star)$ for all $F \in \cL$. Of course, $M_F$ is the localization 
of $M$ at $\{ \, \bx^\bc \mid \bc \in F \cap C \, \}$. We may assume that 
$F = \{ \, u \in \R^d \mid h_1(u) = h_2(u) = \cdots 
= h_n(u) = 0 \, \}$. 

Let $\ba \in \Psi_C$, and take 
$0 \ne \bc \in C$ with $\supp(-\ba) \supset \supp(\bc)$.  
If $\bc \in F$, $M_F$ is a module over $R[\bx^{-\bc}]$. 
Thus the multiplication by $\bx^\bc$ gives a bijection 
$(M_F)_{\ba-\bc}  \to (M_F)_\ba$.  So we may assume that $\bc \not \in F$. 
If $[M_F]_\ba \ne 0$, then we have $\ba = \bb - \ff$ 
for some $\bb \in \Cb$ and $\ff \in F \cap \Z^d$. Thus 
$\supp(-\ba) \subset \supp(\ff) \subset \{ n+1, n+2, \cdots, t \}$. 
On the other hand, since $\bc \not \in  F$, we have $i \in \supp(\bc) 
\subset \supp(-\ba)$ for some $i \leq n$. This is a contradiction. 
Hence $[M_F]_\ba = 0$. Similarly, we can see that  $[M_F]_{\ba -\bc} = 0$.
So any map $[M_F]_{\ba-\bc}  \to [M_F]_\ba$ is bijective.

(Sufficiency) \, By the local duality, the graded Matlis dual 
$H_\m^i(M)^\vee$ of $H_\m^i(M)$ is finitely generated. 
For any $\ba \in \Z^d$ and $0 \ne \bc \in C$, we have 
$- (\ba + n \bc) \in \Phi_C$ for sufficiently large $n$. By the assumption, if 
$[H_\m^i(M)^\vee]_\ba \ne 0$ then $-\ba \not \in \Psi_C$. 
Hence $\m^n H_\m^i(M) = 0$ for sufficiently large $n$, and 
$H^i_\m(M)$ has finite length. 
\qed

\section*{Acknowledgments}
The author is grateful to Professor Ngo Viet Trung for telling him the 
ring given in Example~\ref{square} (2). He also thanks Professors 
Ezra Miller and  Tim R\"omer for careful reading and comments on  
an earlier version of this paper, and Professor Mitsuyasu Hashimoto 
for useful comment around Lemma~\ref{wdel}.


\begin{thebibliography}{99}
\bibitem{BH} W. Bruns and J. Herzog,  Cohen-Macaulay rings, revised edition,  
Cambridge University Press, 1998.

\bibitem{ER}
J.A. Eagon and V. Reiner,  Resolution of Stanley-Reisner rings and Alexander 
duality,  {\it J. Pure and Appl. Algebra} {\bf 130} (1998), 265--275.

\bibitem{D} A.M. Duval, 
Algebraic shifting and sequentially Cohen--Macaulay simplicial complexes, 
{\it Electron. J. Combin.} {\bf 3} (1996), Research Paper 21.

\bibitem{Ei} 
D. Eisenbud, 
The geometry of syzygies: 
A second course in commutative algebra and algebraic geometry, 
Grad. Texts in Math., Vol. 229, 2005. 

\bibitem{GSW}
S. Goto, N. Suzuki and K. Watanabe,  
On affine semigroup rings. 
{\it Japan. J. Math.}  {\bf 2} (1976), 1--12.

\bibitem{GS}
D. Grayson and M. Stillman, Macaulay2 - A system for computation 
in algebraic geometry and commutative algebra, 1997, 
{\tt http://www.math.uiuc.edu/Macaulay2.} 

\bibitem{Hart} R. Hartshorne,
Complete intersections and connectedness, {\it Amer. J. Math.}
{\bf 84} (1962), 497--508.

\bibitem{HM} D. Helm and E. Miller, 
Algorithms for graded injective resolutions and local cohomology over semigroup rings, 
{\it J. Symbolic Comput.} {\bf 39} (2005) 373--395. 

\bibitem{HTT} J. Herzog, Y. Takayama and N. Terai, 
On the radical of a monomial ideal, {\it Arch. Math.} 
{\bf 85} (2005), 397--408. 

\bibitem{Hi} T. Hibi, Level rings and algebras with straightening laws,  
{\it J. Algebra} {\bf 117} (1988), 343--362.

\bibitem{IR} B. Ichim and T. R\"omer, On toric face rings, 
preprint, math.AC/0605150.  

\bibitem{I} B. Iversen,  Cohomology of sheaves. Springer-Verlag, 1986. 

\bibitem{Mil0} E. Miller, 
The Alexander duality functors and local duality with monomial support, 
{\it J. Algebra.} {\bf 231} (2000), 180--234. 

\bibitem{Mil} 
E. Miller, Cohen-Macaulay quotients of normal affine semigroup 
rings via irreducible resolutions, {\it Math. Res. Lett.} 
{\bf 9} (2002), 117-128.  

\bibitem{MS} E. Miller and B. Sturmfels, 
Combinatorial Commutative Algebra, 
Grad. Texts in Math., Vol. 227,  Springer, 2004. 

\bibitem{R1} T. R\"omer,  Cohen-Macaulayness and squarefree modules, 
{\it Manuscripta Math.,} 104 (2001), 39--48.

\bibitem{St} R. Stanley,  
Combinatorics and commutative algebra, 2nd ed. Birkh\"auser 1996.   

\bibitem{St2} R. Stanley, 
Generalized $h$-vectors, intersection cohomology of toric varieties, 
and related results, in Commutative Algebra and Combinatorics 
(M. Nagata and H. Matsumura, eds.), Advanced Studies in Pure Mathematics 
{\bf 11}, North-Holland, 1987, pp. 187--213.

\bibitem{SV} J. St\"uckrad and W. Vogel,  
``Buchsbaum rings and applications," Springer-Verlag, 1986. 

\bibitem{T} Y. Takayama, 
Combinatorial characterizations of generalized Cohen-Macaulay monomial ideals, 
{\it Bull. Math. Soc. Sci. Math. Roumanie} {\bf 48} (2005) 327--344.

\bibitem{Y1} K. Yanagawa,  
Bass Numbers of Local cohomology modules with supports in monomial ideals, 
{\it Math. Proc. Cambridge Philos. Soc.} {\bf 131} (2001), 45--60. . 

\bibitem{Y2} K. Yanagawa, 
Sheaves on finite posets and modules over normal semigroup rings, 
{\it J. Pure and Appl. Algebra} {\bf 161} (2001), 341--366. 

\bibitem{Y6} K. Yanagawa, 
Stanley-Reisner rings, sheaves, and Poincar\'e-Verdier duality,  
{\it Math. Res. Lett.} {\bf 10} (2003) 635--650. 

\bibitem{Y5} K. Yanagawa, 
Derived category of squarefree modules and local cohomology with monomial 
ideal support, {\it J. Math. Soc. Japan} {\bf 56} (2004) 289--308. 
\end{thebibliography}
\end{document}